\documentclass[reqno,11pt]{amsart}

\usepackage[framemethod=tikz]{mdframed}
\usepackage{lipsum}

\usepackage[latin1]{inputenc}
\usepackage[english]{babel}

\usepackage{amscd}

\usepackage{amsmath}
\usepackage{amsthm}
\usepackage{amssymb}

\usepackage[notref, notcite, final]{showkeys}
\usepackage{color}
\usepackage{graphicx}
\usepackage[hidelinks]{hyperref}
\usepackage{cleveref}
\usepackage{fullpage}
\usepackage{comment}
\usepackage[shortlabels]{enumitem}

\usepackage{dsfont}

\newcommand{\pp}{\partial}
\newcommand{\id}{\mathcal{I}}
 		
\newcommand{\nn}{\mathbb{N}}
\newcommand{\rr}{\mathbb{R}}
\newcommand{\cc}{\mathbb{C}}
\newcommand{\hh}{\mathbb{H}}

\newcommand{\Q}{\mathcal{Q}}
\renewcommand{\Re}{\mathrm{Re}}

\newcommand{\boundOP}{\mathcal{B}}

\newcommand{\closOP}{\mathcal{K}}

\newcommand{\dom}{\operatorname{dom}}

\newcommand{\fcts}{\mathfrak{F}}

\usepackage[autostyle]{csquotes}

\newtheorem{theorem}{Theorem}[section]

\newtheorem{proposition}[theorem]{Proposition}

\newtheorem{method}[theorem]{Method}

\theoremstyle{definition}
\newtheorem{definition}[theorem]{Definition}

\theoremstyle{remark}
\newtheorem{remark}[theorem]{Remark}
\crefname{enumi}{}{}

\title{\bf An introduction to hyperholomorphic spectral theories and  fractional powers of vector operators}

\author[F. Colombo]{Fabrizio Colombo}
\address{(FC)
Politecnico di Milano\\Dipartimento di Matematica\\Via E. Bonardi, 9\\20133
Milano, Italy}
\email{fabrizio.colombo@polimi.it}

\author[J. Gantner]{Jonathan Gantner}
\address{(JG)
(was PhD student at) Politecnico di Milano\\Dipartimento di Matematica\\Via E. Bonardi, 9\\20133
Milano, Italy
} \email{jonathan.gantner@gmx.at}

\author[S. Pinton]{Stefano Pinton}
\address{(SP)
Politecnico di Milano\\Dipartimento di Matematica\\Via E. Bonardi, 9\\20133
Milano, Italy
} \email{stefano.pinton@polimi.it}

\begin{document}
\maketitle

\begin{abstract}
The aim of this paper is to give an overview  of the spectral theories associated with
the notions of holomorphicity in dimension greater than one.
A first natural extension is the theory of several complex variables
 whose Cauchy formula is used to define
 the holomorphic functional calculus for $n$-tuples of
 operators $(A_1,...,A_n)$.
A second way
is to consider hyperholomorphic functions of quaternionic or paravector variables.
In this case, by the Fueter-Sce-Qian mapping theorem, we have two different notions of
hyperholomorphic functions that are called slice hyperholomorphic functions and monogenic functions.
Slice hyperholomorphic functions generate the spectral theory based on the $S$-spectrum
 while monogenic functions induce the spectral theory based on the monogenic spectrum.
There is also an interesting relation between the two hyperholomorphic spectral theories via
 the $F$-functional calculus.
 The two hyperholomorphic spectral theories have different and complementary applications.
 Here we also discuss how to define
 the fractional Fourier's law for nonhomogeneous materials, such definition is
    based on the spectral theory on the $S$-spectrum.
\end{abstract}

\noindent AMS Classification 47A10, 47A60.


\noindent {\em }
\date{today}

\section{Introduction}

The problem to define functions of an operator $A$ or of an $n$-tuple of operators $(A_1,...,A_n)$ is very
important both in mathematics and in physics and has been investigated
with different methods starting from the beginning  of the last century.
The spectral theorem is one of the most important tools to define functions
of normal operators on a Hilbert space
and it is of crucial importance in quantum mechanic as well as the Weyl functional calculus.

\medskip
The theory of holomorphic functions plays a central role in operator theory.
In fact, the Cauchy formula
allows to define the holomorphic functional calculus  (often called Riesz-Dunford functional calculus) in Banach spaces \cite{Bourbaki}, and this calculus
 can be extended to unbounded operators.
For sectorial operators the
  $H^\infty$-functional calculus, introduced by A. McIntosh in \cite{MI86}, turned out to be the most important extension.

\medskip
Holomorphic functions of one complex variable
$f:\Omega \subseteq \mathbb{C} \to \mathbb{C}$ (denoted by $\mathcal{O}(\Omega)$) have the following extensions:

\begin{itemize}
\item[(E1)]
Systems of Cauchy-Riemann equations,
 for functions $f:\Pi\subseteq \mathbb{C}^n \to \mathbb{C}$, give the theory of several complex variables.
\item[(E2)]
Holomorphicity of vector fields is connected with quaternionic-valued functions
or more in general with Clifford algebra-valued functions.
There are two different extensions that are obtain by the Fueter-Sce-Qian theorem, also called the
Fueter-Sce-Qian construction, and gives two different notions of hyperholomorphic functions, see for more details \cite{bookSCE}.
\end{itemize}

Consider functions defined on an open set $U$ in the quaternions $\mathbb{H}$ or in $\mathbb{R}^{n+1}$ for Clifford algebra-valued functions,
then the Fueter-Sce-Qian extension consists of two steps.
\begin{itemize}
\item[]Step (I)  gives the class of slice hyperholomorphic functions (denoted by $SH(U)$),
 these functions are also called slice monogenic for Clifford algebra-valued functions and slice regular in the quaternionic case.
\item[]Step (II)  gives the monogenic functions (denoted by $M(U)$) and Fueter regular functions in the case of the quaternions.
\end{itemize}
Both classes of hyperholomorphic functions have a Cauchy formula that can be used to define
functions of quaternionic operators or of $n$-tuples of
operators that do not necessarily commute.
\begin{itemize}
\item[(S)]
The Cauchy formula of slice hyperholomorphic functions generates the $S$-functional calculus for quaternionic linear operators
or for $n$-tuples of not necessarily commuting operators, this calculus is based on the the notion of
$S$-spectrum.
The spectral theorem for quaternionic operators is also based on the $S$-spectrum.
\item[(M)]
The Cauchy formula of monogenic functions generates the monogenic functional calculus that is based on the monogenic spectrum.
\end{itemize}
The hyperholomorphic functional calculi coincide with the  Riesz-Dunford functional calculus
when they are applied to a single complex operator.

\medskip
If we denote by the symbol ${FSQ}$ the Fueter-Sce-Qian construction
then the following diagram illustrates the possible extensions:
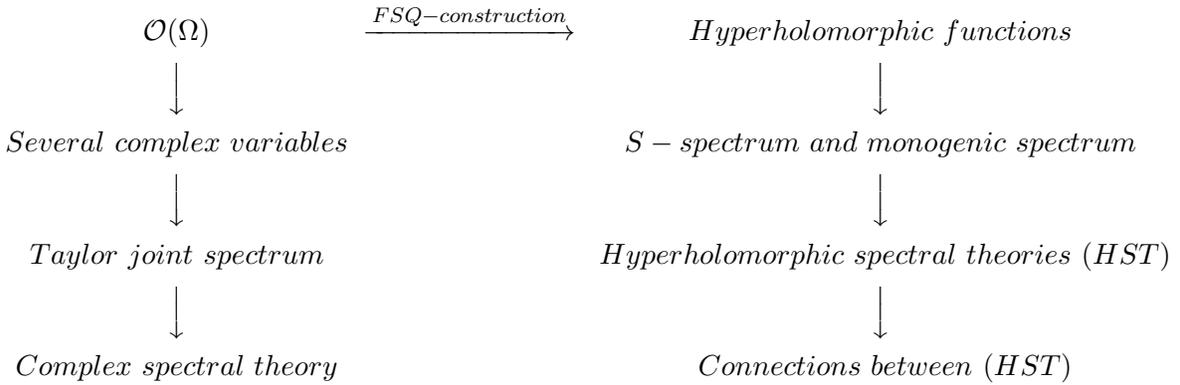
\begin{equation*}
\begin{CD}
\textcolor{black}{\mathcal{O}(\Omega)} @>{FSQ}-construction >>  \textcolor{black}{Hyperholomorphic \ functions} \\   @V VV
  @V VV
\\
\textcolor{black}{Several\ complex \ variables}  @. \textcolor{black}{S-spectrum \ and \ monogenic \ spectrum}
\\
@V VV    @V VV
\\
\textcolor{black}{Taylor \ joint\ spectrum} @. \ \textcolor{black}{ Hyperholomorphic \ spectral \ theories\ (HST)}
\\
@V VV    @V VV
\\
\textcolor{black}{Complex\ spectral\ theory} @. \ \textcolor{black}{Connections\ between\ (HST) }
\end{CD}
\end{equation*}

\medskip
The first mathematicians who understood the importance of hypercomplex analysis to define functions of noncommuting operators on Banach spaces were
A. McIntosh and his collaborators, staring from preliminary results in \cite{KISIL}.
 Using the theory of monogenic functions they developed the monogenic functional calculus and several of its applications, see \cite{6jefferies}.
 The $S$-functional calculus, and in general the spectral theory on the $S$-spectrum, started its development only in 2006
 when F. Colombo and I. Sabadini discovered the $S$-spectrum.
 The discovery of the $S$-spectrum and of the $S$-functional calculus is well explained in the introduction of the book \cite{6CKG} with a complete list of the references and it is also described how hypercomplex analysis methods were used to identify the appropriate notion of  quaternionic
  spectrum  whose existence was suggested by quaternionic quantum mechanics.

  \medskip
 If we denote by $\mathcal{B}(V)$ the Banach space of all bounded right linear operators acting on a two sided quaternionic Banach space $V$
 then the appropriate notion quaternionic spectrum, that is called the $S$-spectrum,
  is defined in a very counterintuitive way because it involves the square of the quaternionic linear operator $T$ and it is define as:
$$
\sigma_S(T)=\{ s\in \mathbb{H}\ | \ T^2-2s_0T+|s|^2\mathcal{I}\ \ {\rm is \ not\ invertible\ in}\ \  B(V)\}.
$$
The $S$-spectrum for quaternionic operators can be naturally defined also for paravector operators when we work in a Clifford algebra,
see \cite{CSSFUNCANAL} and the book \cite{6css}.

\medskip
The problem of the definition of the quaternionic
spectrum for the quaternionic spectral theorem has been an open problem for long time
even though several attempts have been done by several authors in the past decades,
 see e.g. \cite{Teichmueller,Viswanath},
 however the correct definition of spectrum  was not specified.
Finally in \cite{6SpecThm1} the spectral theorem on the $S$-spectrum was proved
 for both bounded and unbounded normal operators on a quaternionic Hilbert space.

\medskip
The main problems with the quaternionic notion of spectrum can be easily described with the following considerations related to
 bounded linear operators just for the sake of simplicity.
Let $T:V\to V$ be a right linear bounded quaternionic operator
  acting on a two sided quaternionic Banach space $V$.
If we readapt the notion of spectrum for a complex linear operator to the quaternionic setting
we obtain two different notions of spectra because of the noncommutativity of the quaternions.
The left spectrum $\sigma_L(T)$ of $T$ is defined as
$$
\sigma_L(T)=\{ s\in \mathbb{H}\ \ | \ \ s\mathcal{I}-T\ \ \ {\rm is\ not\  invertible\ in\ }\mathcal{B}(V) \},
$$
where the notation $s\mathcal{I}$  in $\mathcal{B}(V)$ means that
$(s\mathcal{I} )(v)=sv$.
The right spectrum $\sigma_R(T)$ of $T$ is associated with the right eigenvalue problem, i.e.,
the search of those quaternions $s$ such that there exists a nonzero vector $v\in V$
satisfying
$$
T(v)=vs.
$$
In both spectral problems it is unclear how to associate to the spectrum
a resolvent operator with the property of being an hyperholomorphic function operator-valued.
In fact, for the  left spectrum $\sigma_L(T)$
it is not clear what notion of hyperholomorphicity is associated to the map $s\to (s\mathcal{I}-T)^{-1}$, for $s\in \mathbb{H}\setminus \sigma_L(T)$
 and
for the right spectrum it is even more weird because the operator  $ \mathcal{I}s -T$ (where $\mathcal{I}s$ means $(\mathcal{I}s )(v)=vs$)
   is not linear, so it is not clear which operator is the candidate to be the resolvent operator.

\begin{remark}
{\rm
One of the main motivations that suggested the existence of the $S$-spectrum
  is the  paper \cite{BF} by G. Birkhoff and J. von Neumann,
   where they showed that quantum mechanics can be formulated also on quaternionic numbers.
     Since that time, several papers and books treated this topic, however it is interesting, and somewhat surprising,
     that an appropriate notion of spectrum for quaternionic linear operators was not present in the literature.
     Moreover, in quaternionic quantum mechanics the right spectrum $\sigma_R(T)$ is the most useful notion of spectrum
 to study the bounded states of a quantum systems.
Before 2006 only in one case the quaternionic spectral theorem was proved specifying the spectrum
 and it is the case of quaternionic normal matrices, see \cite{fpp}, where
the  right spectrum $\sigma_R(T)$ has been used.
}
\end{remark}

\medskip
Now we recall some research directions and applications
of the hyperholomorphic function theories and related spectral theories.

\medskip
The {\em first step of FSQ-construction}  generates {\em slice hyperholomorphic functions} and the {\em spectral theory of the $S$-spectrum}, we have:
 \begin{itemize}
\item
The foundation of the quaternionic spectral theory on the $S$-spectrum are organized in the books \cite{6CG,6CKG}, and for paravector operators see \cite{6css}.
\item
The mathematical tools for quaternionic quantum mechanics is the spectral theorem based on the $S$-spectrum \cite{6SpecThm1,JONAQS}.
\item
Quaternionic evolution operators, Phillips functional calculus, $H^\infty$-functional calculus, see \cite{6CG}.
\item
Quaternionic approximation \cite{6GSBook}.
\item
The characteristic operator functions and applications to linear system theory \cite{6COFBook}.
\item
Quaternionic spectral operators \cite{6JONAME}.
\item
 Quaternionic perturbation theory and invariant subspaces \cite{6CCKS}.
\item
Schur analysis in the slice hyperholomorphic setting \cite{6ACSBOOK}.
\item
The theory of function spaces of slice hyperholomorphic functions \cite{6SpacesBook}.
\item
New classes of fractional diffusion problems based on fractional powers of quaternionic linear operators, see the book \cite{6CG} and the more recent contributions \cite{frac4,frac5,frac1,frac2,frac3}.
\end{itemize}
In the last section of this paper we explain how to treat
fractional diffusion problems using the quaternionic spectral theory on the $S$-spectrum
 and we show some of the recent results on
fractional Fourier's law for nonhomogeneous materials recently obtained.
An example of problems that we can treat is the following.

We warn the reader that in this paper, with an abuse of notations,
 we use the symbol $\underline{x}$ for both the coordinates of a point $(x_1,x_2,\ldots, x_n)\in \mathbb{R}^n$
  or for the vector part of a quaternion or for the imaginary part of a paravector in a Clifford algebra.

Let $\Omega$ be a bounded or an unbounded domain in $\mathbb{R}^3$ and let $\tau>0$ and
denote by $v$ the temperature of the material contained in $\Omega$.
Let $\underline{x}=(x_1,x_2,x_3)\in \Omega$ and consider the evolution problem
\begin{equation}\label{PROBT}
\left\{\begin{split}
&
\partial_t v(\underline{x},t)  +{\rm div}\, T(\underline{x}) v(\underline{x},t) = 0,\ \ \ \ (\underline{x},t)\in \Omega\times (0,\tau]
\\
&
v(\underline{x},0)=f(\underline{x}),\ \ \ \underline{x}\in \Omega
\\
&
v(\underline{x},t)=0,\ \ \underline{x}\in \partial \Omega \ \ \ \ t\in [0,\tau],
\end{split}
\right.
\end{equation}
where $f$ is a given datum
and the heat flux for the nonhomogeneous material contained in $\Omega$, is given by the vector differential operator:
\begin{equation}\label{Tnoncom}
T(\underline{x}) = a(\underline{x})\partial_{x_1} e_1 + b(\underline{x}) \partial_{x_2}e_2 + c(\underline{x})\partial_{x_3}e_3,\ \ \ \underline{x}\in \Omega.
\end{equation}
We determine the conditions on the coefficients
 $a$, $b$, $c:\Omega\to \mathbb{R}$
 under which the operator $T(\underline{x})$ generates the fractional powers $P_{\alpha}(T(\underline{x}))$ of $T(\underline{x})$, for $\alpha\in (0,1)$.
 The vector part of  $P_{\alpha}(T(\underline{x}))$ of $T(\underline{x})$ is defined to be
 the nonlocal Fourier's law associated with $T(\underline{x})$.

\medskip
The {\em second step of FSQ-construction}
  generates {\em Fueter or monogenic functions} and the
{\em spectral theory on the monogenic spectrum}. We highlight some references for the research directions in this area:

\begin{itemize}
\item
Monogenic spectral theory and applications \cite{6jefferies}.
Here one can also find the relations of the monogenic functional calculus
with the Taylor functional calculus and the Weyl functional calculus see also some of the original contributions
 \cite{6MQ,6JM,6JMP,6MP,6qian1}.
\item
    Harmonic analysis in higher dimension, singular integrals and Fourier transform
    see the recent book \cite{BOOKTAO}.
\item
Algebraic Analysis of Dirac systems \cite{6DIRAC2}.
\item
The theory of spinor valued function \cite{6DIRAC3}.
\item
Boundary value problems treated with quaternionic techniques \cite{6GURLYSPROSS}.
    \item
    The extension of Schur analysis in the Fueter setting and related topics
 \cite{6AS1,6ASV1,6ASV2}.
 \item
 Dirac operator on manifolds and spectral theory \cite{bookTF,6DIRAC4}.
  \end{itemize}

\medskip
 We conclude by saying that before the recent works on slice hyperholomorphic functions,
 this function theory was simply seen as an intermediate step in the Fueter-Sce-Qian's construction.
 The literature on hyperholomorphic function theories and related spectral theories is nowadays very large.
 For the function theory of slice hyperholomorphic functions the
 main books are \cite{6SpacesBook,6EntireBook,6css,6GSBook,6GSSbook}, while
 for the spectral theory on the $S$-spectrum we mention the books \cite{6COFBook,6CG,6CKG,6css}.
 For the Fueter and monogenic function theory and related topics
 see the books \cite{6DIRAC1,6DIRAC2,6DIRAC3,6DIRAC4,6DIRAC5,6jefferies,6DIRAC6}.

\section{Spectral theory in the complex setting}\label{sec2}

In this section we discuss what is a functional calculus of a single operator on a Banach space and also for the case of several operators.
When we consider a closed linear operator $A$ with domain $\mathcal{D}(A)\subset X$,
where $X$ is a Banach space, the resolvent set $\rho(A)$ of $A$ is defined as
$$
\rho(A)=\{ \lambda\in \mathbb{C}\ |\ (\lambda \mathcal{I}- A)^{-1}\in B(X)\}
$$
 where $B(X)$ in the space of all bounded linear operators on $X$
and the spectrum of $A$ is the set $\sigma(A)=\mathbb{C}\setminus \rho(A)$ and for $\lambda \in \rho(A)$ the map
$\lambda \to (\lambda \mathcal{I}- A)^{-1}$ is called the resolvent operator.
We start with the following intuitive definition of what is a functional calculus.

\medskip
{\it
A functional calculus for a closed linear operator $A$ on complex Banach space $X$ is a mathematical technique that
allows to construct in a meaningful way an operator $f(A)$ for any function $f$ in a certain
class of functions $\fcts$ defined on sets that contain the spectrum $\sigma(A)$ of $A$.
}

\medskip
The formulation {\em in a meaningful way}  usually means that the functional calculus is compatible with
formally plugging $A$ into the function $f$, whenever this is possible.
 That is, whenever $f(z)$ can be expressed by a formula so that formally replacing $z$ by the operator $A$ yields an expression that is meaningful, then $f(A)$ should correspond to this expression. Some examples should clarify this idea:
\begin{enumerate}[(i)]
\item[(a)] For any polynomial $p(z) = a_nz^n+\ldots+a_1z + a_0\in\fcts$ with $a_{\ell}\in\cc$, the operator $p(A)$ should be given by $p(A) = a_nA^n+\ldots + a_1 A + a_0\id$.
\item[(b)] For any $\lambda\in\rho(A)$ with $R_{\lambda}(z) = (\lambda-z)^{-1}\in \fcts$, the functional calculus is compatible with the resolvent operator at $\lambda$. That is, we have $R_{\lambda}(A) = (\lambda\id - A)^{-1}$.
\item[(c)] For any rational function $r(z) = p(z)/q(z)\in \fcts$ with polynomials $p(z) = a_nz^n+\ldots+a_1z + a_0$ and $q(z) = b_nz^n+\ldots + b_1z + b_0$ with $a_{\ell},b_{\ell}\in\cc$, the operator $r(A)$ should be given by $r(A) = p(A)q(A)^{-1}$, where $p(A) = a_nA^n+\ldots+a_1A + a_0\id$ and $q(A) = b_mA^m+\ldots + b_1A + b_0\id$.
\item[(d)] If $A$ is the infinitesimal generator of a strongly continuous group $U_{A}(t),  t\geq 0$ and $\exp(tz)\in\fcts$, then $\exp(tA) = U_{A}(t)$.
\end{enumerate}
Of course in the case of unbounded operators one has to pay attention to the domain of the operators.
Usually the class $\fcts$ constitutes an algebra, often even a Banach algebra,
and the {\em meaningfulness} of the functional calculus as described above follows from the compatibility of the functional calculus with the algebraic operation. Precisely, a functional calculus usually satisfies several (or all) of the following statements:
\begin{enumerate}[(i)]
\item[(I)] The functional calculus is an algebra homomorphism that is
$(af + g)(T) = af(T) + g(T)$ for all $f,g\in\fcts$ and all $a\in\cc$.
\item[(II)]
For $f$ and $g\in \fcts$ such that $fg\in \fcts$ we expect $(fg)(A)=f(A)g(A)$.
\item[(III)] For $f(z) = 1$, we have $f(A) = 1(A) = \id$.
\item[(VI)] For the identity $f(z) = z$, we have $f(A) = z(A) = A$.
\item[(V)] If $X$ is  a Hilbert space and both $f$ and $\overline{f}(z) = \overline{f(z)}$ belong to $\fcts$, then $\overline{f}(A) = f(A)^*$.
\item[(VI)] If $\fcts$ is normed, then the functional calculus defines a continuous mapping into the space of bounded linear operators $\boundOP(X)$ on $X$, that is $\| f(A)\|_{\boundOP(X)} \leq C \| f\|_{\fcts}$.
\end{enumerate}
We restrict ourselves to the case  $\fcts$ consists of functions that are at least continuous and
to the case that the topology on $\fcts$ is coarses than the topology of locally uniform convergence.
In particular, convergence in $\fcts$ implies pointwise convergence.

For the measurable functional calculus
the main statements of this section hold true,
 but it needs different arguments to show them.

There are two main methods for defining a functional calculus, often both of them can be applied
 in order to construct a specific functional calculus.

\begin{method}\label{method1} One considers a subalgebra $\fcts_0$ of $\fcts$ that is dense in $\fcts$ such that $f(A)$ can be defined easily for any $f\in\fcts_0$ (for instance the set of polynomials or the set of rational functions in $\fcts$). If $f\in \fcts$ is arbitrary, one chooses an approximating sequence $(f_n)_{n\in\nn}$ in $\fcts_0$ for $f$ and defines
$$
f(A) := \lim_{n\to +\infty} f_n(A).
$$
\end{method}
\begin{method}\label{method2}
If any $f\in\fcts$ admits an integral representation
\[
f(z) = \int K(\xi,z) f(\xi)\,d\mu(\xi),
\]
and formally replacing $z$ by $A$ in $K(\xi,z)$ yields a meaningful operator $K(\xi,A)$, then one may define
\[
f(A) := \int K(\xi,A) f(\xi)\,d\mu(\xi).
\]
\end{method}
An example for \cref{method1} is the continuous functional calculus.
With  \cref{method2} we define for example
the Riesz-Dunford-functional calculus or  the Philips functional calculus.

\medskip
There is also another concept behind the notion of functional calculus: {\em the operator $f(A)$ should be defined by letting $f$ act on the spectral values of $A$. In particular, this means that
\begin{equation}\label{FCIntuitionComplex}
Ax = \lambda x\qquad\Longrightarrow \qquad f(A) x = f(\lambda)x,\ \ \ {\rm for}\ \ x\in X.
\end{equation}
}
This idea is usually not to much emphasized in the complex setting when one introduces and explains the concept of a functional calculus
and we shall see here in the sequel the reason way this happens.
However, it is this relation that explains why functional calculi are the fundamental
 techniques for investigating linear operators.
 If it doesn't hold, then a functional calculus does not provide any information about the operator
 even though it generates functions of operators.

 Interestingly enough a deep difference between the theory of complex and the theory of quaternionic linear operators
(or more in general for hyperholomorphic spectral theories)
 is revealed here so that in the latter, the relation \eqref{FCIntuitionComplex} needs to be addressed explicitly.

Let us start our considerations by justifying the importance of the relation \eqref{FCIntuitionComplex}. We therefore recall the easiest result that is shown by an application of a functional calculus.
\begin{theorem}
Let $A\in\cc^{m\times m}$. Then $A$ has an eigenvalue and for any polynomial $p\in\cc[n]$ with $p(A) = 0$, the set of eigenvalues of $A$ is contained in the set of roots of $p$.
\end{theorem}

It is obvious that even for the above, very easy, and fundamental result, the fact that the polynomial functional calculus satisfies the relation \eqref{FCIntuitionComplex} is crucial.

The fact that the relation (\ref{FCIntuitionComplex})
 trivially holds true for any known functional calculus in the complex setting
is shown in the next two results, and this is the reason for which it is not explicitly mentioned.

\begin{theorem}\label{thm1}
Let $\Phi: \fcts\to\boundOP(X)$ be a functional calculus for an operator $A$ defined via method \ref{method1}. If \eqref{FCIntuitionComplex} holds true for any function in $\fcts_0$, then the functional calculus $\Phi$ is compatible with \eqref{FCIntuitionComplex}. This is in particular the case if $\fcts_0$ consists of polynomials or rational functions.
\end{theorem}

\begin{theorem}\label{thm2}
Let $\Phi: \fcts\to\boundOP(X)$ be a functional calculus for an operator $A$ defined via method \ref{method2}. If \eqref{FCIntuitionComplex} holds true for $K(\lambda,\cdot)$ for any $\lambda$, then the functional calculus $\Phi$ is compatible with \eqref{FCIntuitionComplex}.
\end{theorem}

We now recalling that the spectral theorem works as a functional calculus.
In the finite dimensional case when we pick
an $n\times n$ matrix $A = (a_{i,j})$ for $i,j=1,\ldots n$ of complex numbers such that $\overline{a_{i,j}}=a_{j,i}$ for
all $i,j = 1,\ldots, n$. The spectrum $\sigma(A)$ consists of eigenvalues
of $A$, that is, complex numbers $\lambda$ for which the equation $Av = \lambda v$ has a nonzero vector
$v\in \mathbb{C}^n$ as a solution. The hermitian matrix $A$ has a unique decomposition as a finite sum
$$
A =\sum_{j=1}^n \lambda_j E_{\lambda_j,A}
$$
where $ E_{\lambda_j,A}$ is the orthogonal projection onto the eigenspace of the eigenvalue $\lambda_j$.
In the case of bounded selfadjoint (or more in general normal operators) operators $A$ acting in Hilbert space the spectral theorem is the most important tool for the complete description of such operators, in fact we have
$$
 A =\int_{\sigma(A)} \lambda dE_{\lambda;A}
 $$
with respect to a spectral measure $E_{\lambda;A}$ associated with $A$. From the spectral theorem we can define
 $f(A)$  by
$$
 f(A) =\int_{\sigma(A)} f(\lambda) dE_{\lambda;A}
$$
for any continuous (but also bounded Borel measurable) function $f : \sigma(A)\to \mathbb{C}$.
The mapping
$f\to f(A)$ is an algebra homomorphism into the space of bounded linear operators.
The spectral theorem can be generalized to the case of
$n$-tuple of commuting bounded selfadjoint operators $(A_1,\ldots, A_n)$, as
$$
f(A_1,\ldots, A_n) =
\int_{\sigma(A_1,\ldots, A_n)} f(\lambda) dE_{\lambda;A_1,\ldots, A_n}
$$
is valid for the joint spectral measure $E_{\lambda;A_1,\ldots, A_n}$ associated with $A_1,\ldots, A_n$.
The joint spectrum
of $A_1,\ldots, A_n$ in $\mathbb{R}^n$ is the support of $E_{\lambda;A_1,\ldots, A_n}$ and
$f : \sigma(A_1,\ldots, A_n) \to \mathbb{C}$ is any bounded Borel measurable function.
The theorem holds more in general for unbounded normal operators but one has to pay attention to the definitions of commutativity in this case, see the book \cite{Schmuedgen}.

In the case we work in a Banach space the most natural way to define functions of bounded (and also of unbounded) operators is
 the Riesz-Dunford functional calculus
$$
f(A) =\frac{1}{2\pi i}\int_{C}(\lambda \mathcal{I} -A)^{-1}f(\lambda) d\lambda
$$
which holds for all  holomorphic functions $f$ defined in a neighborhood of $\sigma(A)$ in the complex plane.
The simple closed contour $C$ surrounds $\sigma(A)$ and is contained in the domain of the function $f$. There
is a natural generalization of Riesz-Dunford functional calculus for $n$-tuples of bounded operators $A_1,\ldots,A_n$
as
$$
f(A_1,\ldots,A_n) =\frac{1}{(2\pi i)^n}\int_{C_1}\ldots \int_{C_n}
(\lambda_1\mathcal{I}- A_11)^{-1}\ldots (\lambda_n\mathcal{I} - A_n)^{-1}f(\lambda_1,\ldots,\lambda_n)d\lambda
$$
where $d\lambda=d\lambda_1\cdots d\lambda_n$
and  $f$ is any  holomorphic function in a neighborhood of $\sigma(A_1)\times \ldots \times\sigma(A_n)$ in $\mathbb{C}^n$. For each
$j = 1,\ldots, n$ the simple closed contour $C_j$ surrounds $\sigma(A_j)$ and $C_1\times \ldots \times\mathbb{C}_n$ is contained
in the domain of $f$ in $\mathbb{C}^n$.
Also when the operators $A_1,\ldots,A_n$ do not commute with each other, the functional calculus makes
sense with any change in the operator ordering of the function
$$
(\lambda_1,\ldots,\lambda_n)\mapsto  (\lambda_1\mathcal{I}- A_11)^{-1}\ldots (\lambda_n\mathcal{I} - A_n)^{-1} \ \ \ {\rm in}  \ \ \
\mathbb{C}^n\setminus (\sigma(A_1)\times \ldots \sigma(A_n)).
$$
For non commuting operators things are in general more complicated and we will not enter into the details here.

\medskip
The material discussed in this section can be found is several classical books.

\medskip
{\em Some remarks in view of the hyperholomorphic spectral theories.}

(I)
In operator theory on the $S$-spectrum, the statements of \Cref{thm1,thm2}  still hold true but just for a subclass of functions.
The conditions that the functions in the dense subspace $\fcts_0$ resp.
the kernel $K(\lambda,A)$ satisfy \eqref{FCIntuitionComplex} is not true in this setting.
It is neither satisfied by the $S$-resolvent operator, nor by the $F$-resolvent operator,
nor by slice hyperholomorphic rational functions with non-real coefficients.
In particular it is not satisfied, whenever the left-linear structure of the space has a prominent role.

(II)
In general the definitions of hyperholomorphic functional calculi ($S$-functional calculus,
the $F$-functional calculus, the monogenic functional calculus)
the product rule, the composition rule or the spectral mapping theorem do not hold.
One needs additional arguments to show that functional calculi based on the $S$-spectrum satisfy \eqref{FCIntuitionComplex}
at least for a subclass of functions namely, the class of intrinsic functions
and for these functions, the problems mentioned before do not occur.

\section{Spectral theories in the hyperholomorphic setting}

At the beginning of the last century several authors started the study of hyperholomorphic functions
and the most popular class of functions are nowadays called
Fueter (or Cauchy-Fueter) regular functions in the case of the quaternions
 and monogenic functions (or functions in the kernel of the Dirac) for Clifford algebra setting.
The second class of hyperholomorphic functions have been developed more
recently, just at the beginning of this century,
and different definitions are possible even though they are not totally equivalent.

\medskip
In the following we will discuss mainly the implications of the Fueter-Sce-Qian
construction in the Clifford setting, the quaternionic setting is similar and we will
use it for the fractional powers of vector operators in the last section of this paper.

\medskip
Let $\rr_n$ be the real Clifford algebra over $n$ imaginary units $e_1,\ldots ,e_n$
satisfying the relations $e_\ell e_m+e_me_\ell=0$,\  $\ell\not= m$, $e_\ell^2=-1.$
 An element in the Clifford algebra will be denoted by $\sum_A e_Ax_A$ where
$A=\{ \ell_1\ldots \ell_r\}\in \mathcal{P}\{1,2,\ldots, n\},\ \  \ell_1<\ldots <\ell_r$
 is a multi-index
and $e_A=e_{\ell_1} e_{\ell_2}\ldots e_{\ell_r}$, $e_\emptyset =1$.
A point $(x_0,x_1,\ldots,x_n)\in \mathbb{R}^{n+1}$  will be identified with the element
$
 x=x_0+\underline{x}=x_0+ \sum_{j=1}^n x_j e_j\in\mathbb{R}_n
$
called paravector and the real part $x_0$ of $x$ will also be denoted by $\Re(x)$.
The imaginary part of $x$ is defined by
 ${\rm Im}(x)=x_1e_1+\ldots+x_ne_n$ and for the sake of simplicity we also use the notation $\underline{x}$ for ${\rm Im}(x)$.
 The conjugate of $x$ is denoted by $\overline{x}=x_0-{\rm Im}(x)$
and the Euclidean modulus of $x$ is given by $|x|^2=x_0^2+\ldots +x_n^2$.
The sphere of purely imaginary paravectors with modulus $1$, is defined by
$$
\mathbb{S}=\{x=e_1x_1+\ldots +e_nx_n\ |\  x_1^2+\ldots +x_n^2=1\}.
$$
The element  $I\in \mathbb{S}$ are such that $I^2=-1$ so we will denote the complex place with imaginary unit $I$ by $\mathbb{C}_I$.
 For this reason the elements of $\mathbb{S}$ are also called
imaginary units.
Given a non-real paravector $x=x_0+{\rm Im} (x)=x_0+J_x |{\rm Im} (x)|$, $J_x:={\rm Im} (x)/|{\rm Im} (x)|\in\mathbb{S}$, we can associate to it the sphere defined by
$$
[x]=\{x_0+J  |{\rm Im} (x)| \ \ | \ \ J \in\mathbb{S}\}.
$$
The set of quaternions will be denoted by $\mathbb{H}$ and the above definitions adapts in this setting in a natural way.
\begin{definition}
Let $U \subseteq \mathbb R^{n+1}$ (or $U \subseteq \mathbb H$ ). We say that $U$ is
\textnormal{axially symmetric}  if, for every $u+Iv \in U$,  all the elements $u+Jv$ for $J\in\mathbb{S}$ are contained in $U$.
\end{definition}
For operator theory the most appropriate definition of
 slice hyperholomorphic functions is the one that comes from the Fueter-Sce-Qian mapping theorem
 because it allows to define functions on axially symmetric open sets.
\begin{definition}\label{CAUSLICE}
 Let $U\subseteq\mathbb{R}^{n+1}$ be an axially symmetric open set
 and let $\mathcal{U}\subseteq\mathbb{R}\times \mathbb{R}$ be such that $x=u+J v\in U$ for all $(u,v)\in\mathcal{U}$.
We say that a function $f:U \to \mathbb{R}_n$ of the form
$$
f(x)=f_0(u,v)+J f_1(u,v)
$$
is left slice hyperholomorphic if
 $f_0$, $f_1$ are $\mathbb{R}_n$-valued differentiable functions such that
 $$
 f_0(u,v)=f_0(u,-v), \ \ \ f_1(u,v)=-f_1(u,-v)\ \ \ {\rm for \ all}\ \  (u,v)\in \mathcal{U}
 $$
 and if $f_0$ and $f_1$ satisfy the Cauchy-Riemann system
 $$
\partial_u f_0-\partial_vf_1=0,\ \ \ \ \
\partial_vf_0+\partial_u f_1=0.
$$
\end{definition}
The above definition adapt naturally to the quaternionic setting.
Since we will restrict just to left slice hyperholomorphic function
on $U$ we introduce  the symbol ${SH}_L(U)$ to denote them.
The subset of intrinsic functions consist of those slice hyperholomophic functions such that
$f_0$, $f_1$ are real-valued and is denoted by $N(U)$.
We recall that right slice hyperholomorphic  functions are of the form
$$
f(x)=f_0(u,v)+f_1(u,v)J
$$
where $f_0$, $f_1$ satisfy the above conditions.
\begin{definition}[Monogenic functions]
Let  $f : U \to \mathbb{R}_n$ be a continuously  differentiable function defined on an open subset $U \subseteq\mathbb{R}^{n+1}$.
 We say that $f$ is (left) monogenic on $U$, if
 $$
 Df(x)=0
 $$
 where $D$ is the Dirac operator defined by
 $$
 D=\partial_{x_0}+\sum_{j=1}^n e_j\partial_{x_j}.
 $$
\end{definition}

The definition of slice hyperholomorphic functions and of monogenic functions
can be seen as to two steps in the Fueter-Sce- Qian constructions to extend holomorphic functions to
 dimension greater than one for  the vector-valued functions (quaternionic or Clifford valued-functions).

In fact, starting from holomorphic functions, R. Fueter  in 1935, see \cite{6fueter}, showed
 an interesting way to generate Cauchy-Fueter regular functions.
More then 20 years later in 1957  M. Sce, see \cite{6sce2},  extended this result in a very pioneering and general way that includes Clifford algebras,
see the English translation of his works in hypercomplex analysis with commentaries collected in the recent book \cite{bookSCE}.

\medskip
In the original construction of R. Fueter the holomorphic
  functions are defined on open sets of the upper half complex plane. This condition can be relaxed
by taking function
$$
g(z) = g_0(u,v)+ig_1(u,v),\ \ \ z=x+iy
$$
 defined in a set $D\subseteq\mathbb C$, symmetric with respect to the real axis such that
 $$
g_0(u,-v)= g_0(u,v) \ \ \ {\rm and} \ \ \  g_1(u,-v)= -g_1(u,v)
$$
 namely if $g_0$ and $g_1$ are, respectively, even and odd functions in the variable $v$.
 Additionally the pair $(g_0,g_1)$ satisfies the Cauchy-Riemann system.
 The above remark holds also for M. Sce's theorem that we state in the following for Clifford algebras.
\begin{theorem}[Sce \cite{6sce2}]\label{SCETH}
Consider the Euclidean space $\rr^{n+1}$ whose elements are identified with
 paravectors $x=x_0+\underline{x}$.
Let $\tilde{f}(z) = f_0(u,v)+if_1(u,v)$ be a  holomorphic function
defined in a domain (open and connected) $D$ in the upper-half complex plane and let
$$
\Omega _D= \{x =x_0+\underline{x}\ \ |\ \  (x_0, |\underline{x}|) \in D\}
$$
be the open set induced by $D$ in $\mathbb{R}^{n+1}$.
The following map
 $$
f(x)=T_{FS1}(\tilde{f}):=\textcolor{black}{f_0(x_0,|\underline{x}|)+\frac{\underline{x}}{|\underline{x}|}f_1(x_0,|\underline{x}|)}
$$
takes the holomorphic functions $\tilde{f}(z)$ and induces the Clifford-valued function $f(x)$.
Then the function
$$
\breve{f}(x):=\textcolor{black}{T_{FS2}} \Big(\textcolor{black}{f_0(x_0,|\underline{x}|)+\frac{\underline{x}}{|\underline{x}|}f_1(x_0,|\underline{x}|)}\Big),
$$
where $T_{FS2}:=\Delta_{n+1}^{\frac{n-1}{2}}$ and $\Delta_{n+1}$ is the laplacian in $n+1$ dimensions,
is in the kernel of the Dirac operator, i.e.,
$$
D\breve{f}(x)=0 \ \ \ {\rm on} \ \ \Omega_D.
$$
\end{theorem}
The case in which the operator $\Delta_{n+1}^{\frac{n-1}{2}}$ has a fractional index has been treated by T. Qian in \cite{6qian}.
Observe that for the Fueter's theorem  the operator $T_{FS2}$ is equal to the laplacian $\Delta$ in $4$ dimensions.
Further developments can be found in \cite{6DIX1,6DIX2,6DIX4} see also the survey \cite{6qianSR}.
\medskip
We can summarize the Fueter-Sce contractions as follows.
Denoting by $\mathcal{O}(D)$ the set of holomorphic functions on $D$, by ${N(\Omega_D)}$ the set
of induced functions on $\Omega_D$ (which turn out to be intrinsic slice hyperholomorphic functions) and by $AM(\Omega_D)$ the set
 of axially monogenic functions on $\Omega_D$ the Fueter-Sce construction can be visualized by the diagram:
\begin{equation*}
\begin{CD}
\textcolor{black}{\mathcal{O}(D)}  @>T_{FS1}>> \textcolor{black}{N(\Omega_D)}  @>\ \   T_{FS2}=\Delta^{(n-1)/2)} >>\textcolor{black}{AM(\Omega_D)},
\end{CD}
\end{equation*}
where $T_{FS1}$ denotes the first linear operator of the Fueter-Sce construction and $T_{FS2}$ the second one.
The Fueter-Sce mapping theorem induces two spectral theories according to the two classes of hyperholomorphic functions it generate.

\medskip
Recently also the problem of construction the inversion of the maps that appear in the Fueter-Sce-Qian extension has been treated, we mention
the papers \cite{6BKQS,6BKQS1,6CoSaSo1,6csso1,6cosaso2,6DIX3},
while a different method to connect slise monogenic and monogenic functions is via the
Radon and dual Radon transform, see \cite{6Radon}.

\begin{remark}
{\rm
The theory of slice hyperholomorphic functions was somewhat abandoned until 2006 when G. Gentili and D. C. Struppa
(inspired by C. G. Cullen \cite{6cullen}) introduced
in \cite{6gentilistruppa} the notion of slice regular functions for the quaternions.
Further developments of the theory of slice regular functions were discussed also in \cite{6adv2} and
the above definition was extended by F. Colombo, I. Sabadini and D.C. Struppa, in \cite{6cssisrael},
(see also \cite{6ISRAEL2,6crelle, 6STRUC})  to the Clifford algebra setting.
Slice regular functions as defined in \cite{6gentilistruppa} and their generalization to the Clifford algebra as in \cite{6cssisrael},
called slice monogenic functions,
 possess good properties on specific open sets that are called axially symmetric slice domains.
  When it is not necessary to distinguish between the quaternionic case and the Clifford algebra case we call these functions slice hyperholomorphic.
 The extension slice hyperholomorphic functions on real alternative algebras can be found in \cite{6ghiloniperotti}.
}
\end{remark}
 \begin{remark}{\rm
 It is also possible to define slice hyperholomorphic functions, as functions in the kernel
  of the first order linear differential operator (introduced in \cite{6Global})
 $$
Gf= \Big(|\underline{x}|^2\frac{\partial }{\partial x_0} +  \underline{x} \sum_{j=1}^n  x_j\frac{\partial }{\partial x_j}\Big)f=0,
$$
where $\underline{x}=x_1e_1+\ldots +x_ne_n$ .
While, a forth way to introduce slice hyperholomorphicity, done in 1998 by G. Laville and I. Ramadanoff in the paper \cite{6LAVRAM}, is
inspired by the Fueter-Sce-Qian mapping theorem. They introduce the
so called Holomorphic Cliffordian functions
defined by the differential equation
$D\Delta^m f = 0$ over $\mathbb{R}^{2m+1}$, where $D$ is the Dirac operator.
Observe that the definition via the global operator $G$ requires less regularity of the functions  with respect to the
 definition in \cite{6LAVRAM}.
}
\end{remark}

\medskip
We now recall the hyperholomorphic Cauchy formulas that are the heart of the hyperholomorphic spectral theories.
It is important to remark that the hypotheses of the following Cauchy formula
are related to the Definition \ref{CAUSLICE} of slice hyperholomorphic functions.

\begin{theorem}[Cauchy formula for slice hyperholomorphic functions]\label{Cauchynuovo}
Let $U \subseteq \mathbb{R}^{n+1}$ be an axially symmetric open set  such that
$\partial (U \cap \mathbb{C}_I)$ is union of a finite number of
continuously differentiable Jordan curves, for every $I\in\mathbb{S}$. Let $f$ be an  $\mathbb R_n$-valued
 slice hyperholomorphic function on an open set containing $\overline{U}$ and, for any $I\in \mathbb{S}$,  we set  $ds_I=-Ids$.
Then, for every $x\in U$, we have:
\begin{equation}\label{integral}
 f(x)=\frac{1}{2 \pi}\int_{\partial (U \cap \mathbb{C}_I)} S_L^{-1}(s,x)ds_I f(s),
\end{equation}
where the slice hyperholomorphic  Cauchy  kernel
is given by
\begin{equation}\label{CACHYKER}
S_L^{-1}(s,x)=-(x^2-2{\rm Re}(s)x+|s|^2)^{-1}
 (x-\overline{s})
\end{equation}
and the value of the integral (\ref{integral}) depends neither on $U$ nor on the  imaginary unit
$I\in\mathbb{S}$.
\end{theorem}
\begin{remark}\label{reminop}{\rm
In the paper \cite{6gentilistruppa} was introduced the notion of slice regularity,
besides the definition, the authors treated power series centered at the origin and some consequences.
Without any tools
the Cauchy formula with slice hyperholomorphic kernel and the representation formula
were originally determine with the following elementary considerations.
To determine the slice hyperholomorphic Cauchy kernel $S_L^{-1}(s,q)$
we observe that from the definition  of slice regularity
 its expansion
$$
S_L^{-1}(s,q):=\sum_{m=0}^\infty q^ms^{-1-m}, \ \ |q|<|s|
$$
is true  when  $q$ and $s$ belong to the same complex plane $\mathbb{C}_I$, for $I\in \mathbb{S}$.
Then we ask ourself what is the closed form of the series in the case $q$ and $s$
do not belong to the same complex plane $\mathbb{C}_I$ observing that
\begin{equation}\label{efdsas}
 \Big(\sum_{m=0}^\infty q^ms^{-1-m}\Big)s-q \Big(\sum_{m=0}^\infty q^ms^{-1-m}\Big)=1
\end{equation}
is true also when $s$ and $q$ do not belong to the same complex plane  $\mathbb{C}_I$.
In the quaternionic case it was observed that the inverse $S$ of $S_L^{-1}(s,q)$
 is the non trivial solution of the quaternionic equation
$$
S^2+Sq-sS=0,
$$
which easily follows from (\ref{efdsas}).
The unknown $S$ was determined using the Niven's Algorithm as it is shown in the historical Note 4.18.3 in the book \cite{6css} and it gives
$$
S(s,q)=(q-\overline{s})^{-1}s(q-\overline{s})-q.
$$
Taking the inverse of $S(s,q)$ we have the Cauchy kernel defined in (\ref{CACHYKER}) for the quaternions.
 This strategy to determine the Cauchy kernel shows that in the Clifford setting the Cauchy kernel remains the same
  if we consider
 $$
S_L^{-1}(s,x):=\sum_{m=0}^\infty x^ms^{-1-m}, \ \ |x|<|s|
$$
where $x=x_0+x_1e_1+...+x_ne_n$ and $s=x_0+s_1e_1+...+s_ne_n$ are paravectors.
Moreover, observe that a direct computation of the integral
$$
\frac{1}{2 \pi}\int_{\partial (U \cap \mathbb{C}_I)} S_L^{-1}(s,x)ds_I f(s)
$$
by computing the residues of the singularities of the kernel $S_L^{-1}(s,x)$ in the complex plane $\mathbb{C}_I$,
gives:
\[
\begin{split}
 \frac{1}{2 \pi}\int_{\partial (U \cap \mathbb{C}_I)} S_L^{-1}(s,x)ds_I f(s)
=\frac{1}{2}\Big[   f(u+Jv)+f(u-Jv)\Big] +I\frac{1}{2}\Big[ J[f(u-Jv)-f(u+Jv)]\Big],
\end{split}
\]
choosing any $J\in \mathbb{S}$ for all $x=u+Iv \in U $.
From here one can clearly see the existence of the structure formula (or representation formula) for slice monogenic functions,
see for example \cite{6STRUC} (or \cite{CAUMON}):
$$
f(u+Iv)=\frac{1}{2}\Big[   f(u+Jv)+f(u-Jv)\Big] +I\frac{1}{2}\Big[ J[f(u-Jv)-f(u+Jv)]\Big].
$$
The quaternionic setting is just a particular case and from these observations started a full development
 of the theory of slice hyperholomorphic functions.
}
\end{remark}

The second ingredient for our discussion in the following is the
Cauchy formula for monogenic functions.
\begin{theorem}[Cauchy formula for monogenic functions]
Let $U \subset \mathbb{R}^{n+1}$ be an open set with smooth boundary $\partial U$ and let
 $\eta(\omega)$ be the outer unit normal to $\partial U$ and $dS(\omega)$ be the scalar element of surface area on $\partial U$.
Let $f$ be a monogenic function on an open set that contains $\overline{U}$ then
$$
f(x)=\int_{\partial U} G_\omega(x)\eta(\omega) f(\omega)dS(\omega)
$$
 for every $x$ in $U$,
where the monogenic Cauchy kernel is given by
 $$
G_\omega(x):=\frac{1}{\sigma_n}\frac{\overline{\omega-x}}{|\omega-x|^{n+1}},\ \ \ x,\  \omega\in \mathbb{R}^{n+1}, \ \ x\not=\omega
$$
and
$
\sigma_n:= 2\pi^{\frac{n+1}{2}}/ \Gamma\Big(\frac{n+1}{2}\Big)
$
is the volume of unit n-sphere in $\mathbb{R}^{n+1}$.
\end{theorem}

\medskip
Before to introduce the basic fact on the hyperholomorphic spectral theories we need some important considerations.

\medskip
(I) Holomorphic functions of one complex variable and harmonic analysis are strongly connected since the Cauchy-Riemann operator factorizes the Laplace operator.
 The holomorphic functional calculus and the spectral theorem are based on the same notion of spectrum.

\medskip
(II)  In order to restore the analogy with the holomorphic functional calculus and the spectral theorem in the quaternionic setting we have to replace the classical spectrum with the $S$-spectrum. In, fact the $S$-functional calculus and the quaternionic spectral theorem are both based on the $S$-spectrum.
The Dirac operator factorizes the Laplace operator
the monogenic functional calculus,
 based on the monogenic spectrum, has applications in harmonic analysis and in other related fields.

\medskip
Let us consider a Banach space $V$ over
$\mathbb{R}$  with norm $\|\cdot \|$.
It is possible to endow $V$ with an operation of multiplication by
elements of $\mathbb{R}_n$ which gives a two-sided module over
$\mathbb{R}_n$ and by $V_n$ we indicate the two-sided Banach module  over $\mathbb{R}_n$ given by $V\otimes \mathbb{R}_n$.

We start with the definition of a functional calculus for $(n+1)$-tuples of
not necessarily commuting operators using slice hyperholomorphic functions.
So we consider the paravector  operator
$$
T=T_0+\sum_{j=1}^ne_jT_j,
$$
where $T_\mu\in\mathcal{B}(V)$ for   $\mu=0,1,...,n$,
and where $\mathcal{B}(V)$ is the space of all bounded $\mathbb{R}$-linear operators acting on $V$.
The notion of $S$-spectrum
follows from the Cauchy formula of slice hyperholomorphic functions
and from some not trivial considerations
on the fact that we can replace in the Cauchy kernel $S^{-1}_L(s,x)$ the paravector $x$ by the paravector operator $T$ also in the case the
components $(T_0, T_1,...,T_n)$ of $T$ do not commute among themselves.

\begin{remark}{\rm
We make a crucial observation which justifies the definition of $S$-spectrum and of $S$-resolvent operator.
With the procedure of Remark \ref{reminop}  it is natural to replace the paravector $x$ by the paravector operator
  $T=T_0+T_1e_1+...+T_ne_n$ with bounded not necessarily commuting components $T_\ell$, $\ell=0,...,n$ in the Cauchy kernel series. We obtain
$$
\sum_{m=0}^\infty T^ms^{-1-m}=-(T^2-2{\rm Re}(s)T+|s|^2\mathcal{I})^{-1} (T-\overline{s}\mathcal{I}), \ \ \|T\|<|s|.
$$
even though the components of $T$ do not commute.
From this relation we justify the definition of the $S$-resolvent operator and of the $S$-spectrum.
The quaternionic setting is just a particular case.
}
\end{remark}

We have
 the following  definition.
\begin{definition}[$S$-spectrum]\label{defspscandres}
\index{S-spectrum}\index{S-resolvent set}
Let $T\in\mathcal{B}(V_n)$ be a paravector operator.
We define the $S$-spectrum $\sigma_S(T)$ of $T$  as:
$$
\sigma_S(T)=\{ s\in \mathbb{R}^{n+1}\ \ :\ \ T^2-2 \, {\rm Re}\, (s)T+|s|^2{\mathcal{I}}\ \ \
{\rm is\ not\  invertible\ in \ }\mathcal{B}(V_n)\}
$$
where $\mathcal{I}$ denotes the identity operator.
The $S$-resolvent set of $T$ is defined as
$$
\rho_S(T)=\mathbb{H}\setminus\sigma_S(T).
$$
\end{definition}

\begin{definition} Let  $T\in\mathcal{B}(V_n)$ be a paravector operator and $s\in \rho_S(T)$.
We define the left $S$-resolvent operator as \begin{equation}\label{quatSresolrddlft}
S_L^{-1}(s,T):=-(T^2-2{\rm Re}\,( s) T+|s|^2{\mathcal{I}})^{-1}(T-\overline{s}{\mathcal{I}}).
\end{equation}
\end{definition}
A similar definition can be given for the right resolvent operator.
\begin{definition}\label{ggggg}
We denote by
${SH}^L_{\sigma_S(T)}$
the set of slice hyperholomorphic functions
defined on the axially symmetric set $U$ that contains the $S$-spectrum of $T$.
\end{definition}
 A crucial result for the definition of the $S$-functional calculus is that integral
\begin{equation}\label{quatinteg311}
{\frac{1}{2\pi }} \int_{\partial (U\cap \mathbb{C}_{I})} S_L^{-1} (s,T)\  ds_{I} \ f(s),\ \ \ \ \   {\rm for}\ \ \ \ \ f\in {SH}^L_{\sigma_S(T)}
\end{equation}
 depends neither on $U$ nor on the  imaginary unit
$I\in\mathbb{S}$, so the  $S$-functional calculus turns out to be well defined.

\begin{definition}[$S$-functional calculus]\label{quatfunccalleftright} Let $T\in{B}({V}_n)$ and let $U\subset \mathbb H$ be as above.
   We set $ds_I=- I ds$ and we define the $S$-functional calculus as
\begin{equation}\label{quatinteg311def}
f(T):={{1}\over{2\pi }} \int_{\partial (U\cap \mathbb{C}_I)} S_L^{-1} (s,T)\  ds_I \ f(s), \ \ \ \ \   {\rm for}\ \ \ \ \ f\in {SH}^L_{\sigma_S(T)}.
\end{equation}
\end{definition}
Observe that the definition of the $S$-functional calculus is very natural for non commuting operators in
noncommutative spectral theory. The heart of the general version of the $S$-functional calculus can be found in the original papers \cite{6newresol,CAUMON,JGA,CLOSED} and its commutative version \cite{6FUNC3}.

\medskip
{\em Warning}.
In the monogenic setting the natural functional calculus is for vector operators that is
when we set $T_0=0$ in the paravector operator $T=T_0+\sum_{j=1}^ne_jT_j$. The reason will be clear in the sequel, but to point out this fact we use the symbol
$A = (A_1,\ldots,A_n)$ or $A=\sum_{j=1}^ne_jA_j$ instead of  $ (T_1,\ldots,T_n)$ or $T=\sum_{j=1}^ne_jT_j$.

\medskip
Using the Cauchy integral formula  for monogenic functions, we establish the monogenic functional calculus
 for the $n$-tuple $A = (A_1,\ldots,A_n)$ of bounded linear operators on a Banach space
$X$ by substituting the $n$-tuple $A$ for the vector $x\in \mathbb{R}^n$.

In the following for the monogenic functional calculus we limit ourselves to the most simple case when
 $n$ is odd and the $n$-tuple $A = (A_1,\ldots,A_n)$ of bounded linear operators commute among themselves.
Such restrictions can be removed but one needs to do further considerations.

\begin{remark}
If $n$ is odd, A is a commutative
$n$-tuple, that is, $A_jA_k = A_kA_j$ for $j, k = 1,\cdots, n$, and each operator $A_j$ has real spectrum
$\sigma(A_j) \subset \mathbb{R}$ for $j = 1,\cdots, n$, then for suitable $\omega\in  \mathbb{R}^{n+1}$, the expression
\begin{equation}\label{ggg}
G_\omega(A):=\frac{1}{\sigma_n}\frac{\overline{\omega \mathcal{I}-A}}{|\omega \mathcal{I}-A|^{n+1}}
\end{equation}
makes sense as an element of $\mathcal{B}(V_n)$ and it is called the monogenic resolvent.
\end{remark}
\begin{remark}{\rm
 For an even integer $m$ we have
$$
|\omega \mathcal{I}-A|^{-m}=\Big(\, \Big(\omega_0^2\mathcal{I}+\sum_{j=1}^n(\omega_j\mathcal{I}-A_j)^2\Big)^{-1}\, \Big)^{m/2}
$$
and
$$
\overline{\omega \mathcal{I}-A}=\omega_0I-\sum_{j=1}^n(\omega_j\mathcal{I}-A_j)e_j
$$
for
$\omega=\omega_0+\sum_{j=1}^n\omega_j$.
Observe that the operator
$$
\omega_0^2\mathcal{I}+\sum_{j=1}^n(\omega_j\mathcal{I}-A_j)^2
$$
is invertible in $\mathcal{B}(V)$ for each $\omega_0\not=0$.
}
\end{remark}
\begin{definition}[Monogenic spectrum]
 The function
$$
\omega\mapsto G_\omega(A),
$$
 is  defined on the set $\mathbb{R}^{n+1}\setminus (\{0\}\times \gamma(A))$ where
$$
 \gamma(A)=\{ (\omega_1,\ldots,\omega_n) \ |\ \sum_{j=1}^n(\omega_j\mathcal{I}-A_j)^2\ {\rm is\ not\ invertible\ in}\ \mathcal{B}(V)\}
$$
is called the monogenic spectrum.
\end{definition}
\begin{definition}[The monogenic functional calculus]
Let  $n$  be an odd number and let us assume that
 $A = (A_1,\ldots,A_n)$ is a commutative
$n$-tuple of bounded linear operators (that is $A_jA_k = A_kA_j$ for $j, k = 1,\cdots, n$),
and each operator $A_j$ has real spectrum
$\sigma(A_j) \subset \mathbb{R}$ for $j = 1,\cdots, n$.
If $f$ is a monogenic function on an open set that contains $\overline{U}\subset \mathbb{R}^{n+1}$ with $\gamma(A) \subset U$.
Then we define the monogenic functional calculus as
$$
f(A)=\int_{\partial U} G_\omega(A)\eta(\omega) f(\omega)dS(\omega)
$$
 where $G_\omega(A)$ is the monogenic resolvent operator (\ref{ggg}),  $\eta(\omega)$ is the outer unit normal to $\partial U$ and $dS(\omega)$ is the scalar element of surface area on $\partial U$.
\end{definition}

\medskip
In the case $m = 2, 4, 6,...$ the operator $|\omega \mathcal{I}-A|^{-m}$ needs to be defined in a suitable way.
The direct formulation employs Taylor's functional calculus, but by using the plane wave
decomposition of the Cauchy kernel, the case of even n and noncommuting
operators can be treated simultaneously.
For an $n$-tuple $(A_1,....,A_n)$ of commuting bounded linear operators on a Banach
space $V$ with real spectra, the nonempty compact subset
$ \gamma(A)$ of $\mathbb{R}^n$ coincides with Taylor's
joint spectrum defined  in terms of the Koszul complex.

 The Cauchy formula of slice hyperholomorphic functions allows to define the notion of $S$-spectrum, while
 the Cauchy formula for monogenic functions induces the notion of monogenic spectrum, as illustrated by the diagram:
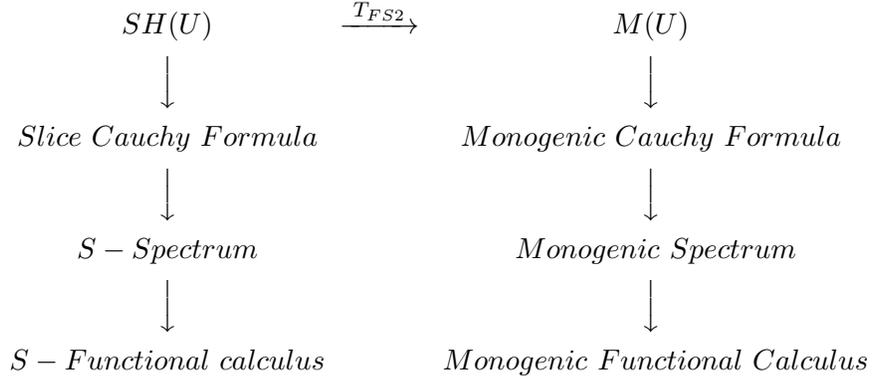
\begin{equation*}
\begin{CD}
\textcolor{black}{SH(U)} @>T_{FS2} >>  \textcolor{black}{M(U)} \\   @V VV
  @V VV
\\
\textcolor{black}{Slice\ Cauchy \ Formula}  @. \textcolor{black}{Monogenic \  Cauchy \ Formula}
\\
@V VV    @V VV
\\
\textcolor{black}{S-Spectrum} @. \ \textcolor{black}{Monogenic \ Spectrum}
\\
@V VV    @V VV
\\
\textcolor{black}{S-Functional \ calculus} @. \ \textcolor{black}{Monogenic\ Functional\ Calculus}
\end{CD}
\end{equation*}
In the above diagram we have replaced the set of  intrinsic functions $N$ by the larger set of slice hyperholomorphic functions $SH$.
This is clearly possible because the map
$T_{FS2}$ is the Laplace operator or its powers.

\medskip
We finally recall that the quaternionic spectral theorem is based on the $S$-spectrum and not on the monogenic spectrum.
In 2015 (and published in 2016)  the quaternionic spectral theorem for quaternionic normal operators was finally proved, see
 \cite{6SpecThm1} (see also \cite{6SpectThmIrene}). Later on  perturbation results of quaternionic normal operators were proved in \cite{6CCKS}.
 Beyond the spectral theorem there are more recent developments in the direction of the characteristic operator functions, see
 \cite{6COFBook} and the theory of quaternionic spectral operators was developed in \cite{6JONAME}.

\medskip
Finally, we wish to give an idea of the structure of the quaternionic spectral theorem. For a complete treatment see \cite{6CKG}.
 If $T \in \mathcal{B}(\mathcal{H})$ is a bounded normal  quaternionic linear operator, on a quaternionic Hilbert space $\mathcal{H}$,
then there exist three quaternionic linear operators $A$, $\mathfrak{J}$, $B$ such that
$T = A + \mathfrak{J} B$, where
$A$ is self-adjoint and $B$ is positive,
 $\mathfrak{J}$ is an anti self-adjoint partial isometry (called imaginary operator). Moreover, $A$, $B$ and $\mathfrak{J}$ mutually commute.
\\
There exists a unique spectral measure $E_I$ on $\sigma_S(T)\cap\mathbb{C}_I^+$ so that for any slice continuous intrinsic function $f  = f_0 + f_1 I$ we have:
\begin{equation}
\langle f(T)x, y \rangle = \int_{\sigma_S(T)\cap\mathbb{C}_I^+} f_0(q) \, d \langle E_I(q)x, y \rangle + \int_{\sigma_S(T)\cap\mathbb{C}_I^+} f_1(q)  \, d \langle \mathfrak{J} E_I(q)x, y \rangle, \quad x,y \in \mathcal{H}.
\end{equation}
This theorem extends to the case of unbounded operators as well and holds true
for a larger class of functions that are not necessarily continuous.

\section{Interaction of the hyperholomorphic spectral theories}

\medskip
Now we formulate the Fueter-Sce-Qian theorem in integral form and we use it
 to define the $F$-functional calculus. This gives a version
of the monogenic functional calculus for $n$-tuples of commuting operators but it is based on the $S$-spectrum instead of the monogenic spectrum.
This calculus was introduced in \cite{6CoSaSo} and further investigated in \cite{6FUNC1,6FUNC2}.

\medskip
It is important to recall that the monogenic functional calculus is defined for $n$-tuples of operators
$A_j$, $j=1,...,n$ that have real spectrum considered as operators $A_j:V\to V$ on the real Banach space $V$.
The $F$-functional calculus has advantages and disadvantages with respect to the monogenic functional calculus.
Precisely,
the $F$-functional calculus allows to consider a much larger class of operators
 because it does not require that the spectrum of the operators $A_j$, $j=1,...,n$  has to be real.
 Moreover, this calculus allows to consider paravector operators
 and not only vector operators as the monogenic functional calculus imposes.

On the other hand,
from the hyperholomorphic functions point of view the $F$-functional calculus
  is less general with respect to the monogenic functional calculus because it works for the subset of monogenic function
 given by
  $$
  \breve{M}=\{\breve{f}\ | \   \breve{f}(x)=\Delta^{\frac{n-1}{2}} f(x) \ {\rm for} \  f\in SH(U)\}.
  $$

\medskip
We now show how the Fueter-Sce mapping theorem provides an
alternative way to define the functional calculus for monogenic functions.
The main idea is to apply the Fueter-Sce operator $T_{FS2}$ to the slice hyperholomorphic Cauchy kernel
 as illustrated by the diagram:
\begin{equation*}
\begin{CD}
{SH(U)} @.  {AM(U)} \\   @V  VV
  @.
\\
{Slice\ Cauchy \ Formula}  @> T_{FS2}>> {Fueter-Sce\ theorem \ in \  integral\  from}
\\
@V VV    @V VV
\\
{S-Functional \ calculus} @. {F-functional \ calculus}
\end{CD}
\end{equation*}
This method generates an integral transform, called the Fueter-Sce mapping theorem in integral form,
 that allows to define the so called $F$-functional calculus.
 This calculus uses slice hyperholomorphic functions and the commutative version of the $S$-spectrum and
 now we show how it works.
 We point out that the operator $T_{FS2}$ has a kernel and one has to pay attention to this fact with the definition of the $F$-functional calculus, more details are given in \cite{6CKG}.

 \medskip
 Now observe that one can apply the powers of the Laplace operators to both sides of \eqref{integral} so that we have
$$
\Delta^h f(x)=\frac{1}{2 \pi}\int_{\partial (U \cap \mathbb{C}_I)} \Delta^h S_L^{-1}(s,x)ds_I f(s).
 $$
 In general, it is not easy to compute $\Delta^h f$ and when we apply $\Delta^h$ to the Cauchy kernel
 written in the form (\ref{CACHYKER}), we do not get a simple formula.
 However, $S_L^{-1}(s,x)$ can be written in two equivalent ways as follows.
\begin{proposition}\label{uguaglianza}
Let $x$, $s\in \mathbb R^{n+1}$ (or in $\mathbb H$ in the quaternionic case)
be such that $x^2 -2x {\rm Re} (s) +|s|^2\not=0$. Then the following identity holds:
\begin{equation}\label{second}
\begin{split}  S_L^{-1}(s,x)&=
-(x^2 -2x {\rm Re} (s) +|s|^2)^{-1}(x-\overline s)
=(s-\bar x)(s^2-2{\rm
Re}(x)s+|x|^2)^{-1}.
\end{split}
\end{equation}
\end{proposition}
If we use the second expression for the Cauchy kernel we find a very simple expression for $\Delta^h S_L^{-1}(s,x)$. In fact, we have:
\begin{theorem}\label{Laplacian_comp}
Let $x$,
$s\in \rr^{n+1}$
be such that
 $x^2 -2x {\rm Re} (s) +|s|^2\not=0$.
Let
$$
S_L^{-1}(s,x)=(s-\bar x)(s^2-2{\rm Re}(x) s+|x|^2)^{-1}
$$
 be the slice monogenic Cauchy kernel and let $\Delta=\sum_{i=0}^n\frac{\partial^2}{\partial x_i^2}$ be the Laplace operator in the variables $(x_0,x_1,...,x_n)$.
Then, for $h\geq 1$, we have:
\begin{equation}\label{hLaplacian}
\Delta^hS_L^{-1}(s,x)=C_{n,h}\ (s-\bar x)(s^2-2{\rm Re}(x)s +|x|^2)^{-(h+1)},
\end{equation}
where
$$
C_{n,h}:=(-1)^h\prod_{\ell=1}^h(2\ell) \prod_{\ell=1}^h (n-(2\ell -1)).
$$
\end{theorem}
The function $\Delta^hS^{-1}(s,x)$ is slice hyperholomorphic in $s$ for any $h\in\mathbb N$ but is monogenic in $x$ if and only if  $h=(n+1)/2$,
 namely if and only if $h$ equals the Sce's exponent. We define the kernel
\[
\begin{split}
\mathcal{F}_L(s,x)&:=\Delta^{\frac{n-1}{2}}S_L^{-1}(s,x)
=\gamma_n(s-\bar x)(s^2-2{\rm Re}(x)s +|x|^2)^{-\frac{n+1}{2}},
\end{split}
\]
where
\begin{equation}\label{gammn}
\gamma_n:=(-1)^{(n-1)/2}2^{(n-1)/2}(n-1)!\Big( \frac{n-1}{2}\Big)!
\end{equation}
which can be used to obtain the Fueter-Sce mapping theorem in integral form.

\begin{theorem}
Let $n$ be an odd number.
Let $f$ be a slice hyperholomorphic function defined in an open set that contains $\overline{U}$, where $U$ is a
bounded  axially symmetric open set. Suppose that the boundary of $U\cap \mathbb{C}_I$  consists of a
finite number of rectifiable Jordan curves for any $I\in\mathbb{S}$.
Then, if $x\in U$, the function $\breve{f}(x)$, given by
$$
\breve{f}(x)=\Delta^{\frac{n-1}{2}}f(x)
$$
is monogenic and it admits the integral representation
\begin{equation}\label{Fueter}
 \breve{f}(x)=\frac{1}{2 \pi}\int_{\pp (U\cap \mathbb{C}_I)} \mathcal{F}_L(s,x)ds_I f(s),\ \ \ ds_I=ds/ I,
\end{equation}
where the integral depends neither on $U$ nor on the  imaginary unit
$I\in\mathbb{S}$.
\end{theorem}
 In the sequel, we will consider bounded paravector operators $T$,
   with commuting components $T_\ell\in\mathcal{B}(V)$ for $\ell=0,1,\ldots ,n$.
Such subset of  ${\mathcal{B}(V_n)}$
will be denoted by $\mathcal{BC}^{\small 0,1}(V_n)$.
\index{$\mathcal{BC}^{\small 0,1}(V_n)$}
The $F$-functional calculus is based on the commutative version of the $S$-spectrum given by
$$
\sigma_S(T)=\{ s\in \mathbb{R}^{n+1}\ \ :\ \ s^2\mathcal{I}-(T+\overline{T})s +T\overline{T}\ \ \
{\rm is\ not\  invertible\ in \ }\mathcal{B}(V_n)\}
$$
where the operator $\overline{T}$ is defined by
$$
\overline{T}=T_0-T_1e_1 - \dots  - T_n e_n.
$$
We observe that for historical reasons the commutative version of the $S$-spectrum is sometimes called $F$-spectrum because it is used for the $F$-functional calculus.
 So we define the $F$-resolvent operators.
\begin{definition}[$F$-resolvent operators]\index{$F$-resolvent operators}
Let $n$ be an odd number and let
$T\in\mathcal{BC}^{\small 0,1}(V_n)$.
 For $s\in \rho_S(T)$
we define the left $F$-resolvent operator by
\begin{equation}\label{FresBOUNDL}
F_L(s,T):=\gamma_n(s{I}-\overline{ T})(s^2\mathcal{I}-(T+\overline{T})s +T\overline{T})^{-\frac{n+1}{2}},
\end{equation}
and
the constants $\gamma_n$ are given in (\ref{gammn}).
\end{definition}

\begin{definition}[The $F$-functional calculus for bounded operators]\label{DEFF_FUNC}
Let $n$ be an odd number, let
 $T = T_0+T_1e_1 + \dots  + T_n e_n\in\mathcal{BC}^{\small 0,1}(V_n)$  and set $ds_I=ds/I$, for $I\in \mathbb{S}$.
 Let ${SH}^L_{\sigma_S(T)}$ and  $U$ be as in Definition \ref{ggggg}.
  We define
\begin{equation}\label{DefFCLUb}
\breve{f}(T):=\frac{1}{2\pi}\int_{\pp(U\cap \mathbb{C}_I)} F_L(s,T) \, ds_I\, f(s).
\end{equation}
\end{definition}
The definition of the $F$-functional calculus is well posed since
the integrals in (\ref{DefFCLUb})  depends neither on $U$ and nor on the imaginary unit $I\in\mathbb{S}$.

\medskip
We conclude this section with some considerations on the hyperholomorphic functional
calculi to show the difference with respect the the complex case.

\medskip
(I)
The product rule holds for the $S$-functional calculus but just in the case  one of the two functions is intrinsic function.
For the monogenic functional calculus the product rule does not hold.
This is due to the fact the that product of two monogenic functions is not monogenic.
For the $F$- functional calculus the product rule does not hold.

\medskip
(II) Regarding the compatibility with polynomials we have:
that the $S$-functional calculus and the monogenic functional calculus are compatible
with slice hyperholomorphic polynomials and with monogenic polynomials, respectively.
For the $F$-functional calculus the compatibility with polynomials holds if we consider
$$
\breve{P}(q)=\Delta P(q)
$$
 where $\breve{P}(q)$ is a monogenic (or Fueter) and $P$ is a slice monogenic polynomials
$$
q\to  T \Rightarrow \breve{P}(q)\Rightarrow \breve{P}(T)
$$

\medskip
(III) The spectral properties of the operator $T$ can be deduced by the $S$-functional calculus and the quaternionic spectral theorem  for which
\begin{equation}
Tx = \lambda x\qquad\Longrightarrow \qquad f(T)x = f(\lambda)x
\end{equation}
when we use intrinsic functions.

\section{The $S$-spectrum approach to fractional diffusion problems}

An important extension of the $S$-functional calculus to unbounded sectorial operators is the
$H^\infty$-functional calculus which is one of the ways to define functions of unbounded operators.
   The  $H^\infty$-functional calculus has been used to define fractional powers of paravector operators and of
   quaternionic linear operators
   that define fractional Fourier laws for nonhomogeneous material in the theory of heat propagation.
   For the original contributions on fractional powers of vector operators and of quaternionic operators and of the $H^\infty$-functional calculus based on the $S$-spectrum see \cite{6hinfty,6MILANO,6Transaction}.
   For a systematic and recent treatment of quaternionic spectral theory on the $S$-spectrum
   and the fractional diffusion problems
   based on techniques on the $S$-spectrum see the books \cite{6CG,6CKG} published in 2019. Moreover, in the monograph
    \cite{6css},   published 2011,  one can find also the foundations of the spectral theory on the $S$-spectrum for $n$-tuples of noncommuting operators.

The theory on the fractional powers of quaternionic operators has been recently applied to physical problems and in
particular to generate the fractional Fourier law for the heat equation that is collected in the papers
\cite{frac1,frac2,frac3,frac4,frac5}, here we give an overview of some of our results.

\medskip
We denote by $\underline x:=(x_1,x_2,x_3)$ a generic point in $\rr^3$ (we warn the reader that the symbol $\underline x$ is also used to denote the imaginary part of a quaternion. It will be clear from the context which is the meaning of the symbol $\underline x$ otherwise it will be specified). Let $\Omega\subset\rr^3$ bounded or unbounded domain (with $\mathcal{C}^1$ boundary), the heat equation for nonhomogeneous materials with the associated initial-boundary conditions descibes the evolution of the heat.
Precisely, we determine $v: \Omega\times (0,\tau]\to \rr$ (for $\tau>0$) such that
\begin{equation}\label{PROBgraddirichlet}
\left\{\begin{split}
&
\partial_t v(\underline{x},t)  +{\rm div}\, T(x) v(\underline{x},t) = 0,\ \ \ \ (\underline{x},t)\in \Omega\times (0,\tau]
\\
&
v(\underline{x},0)=f(\underline{x}),\ \ \ \underline{x}\in \Omega
\\
&
v(\underline{x},t)=0,\ \ \underline{x}\in \partial \Omega \ \ \ \ t\in [0,\tau],
\end{split}
\right.
\end{equation}
where $f$ is a given datum and
\begin{equation}\label{T}
T(\underline x)=\left(
\begin{split}
a_1(\underline x)\partial_{x_1}\\
a_2(\underline x)\partial_{x_2}\\
a_3(\underline x)\partial_{x_3}
\end{split}
\right)
\end{equation}
where we suppose that the coefficients  $a_1$, $a_2$, $a_3: \overline{\Omega} \subset \rr^3\to \mathbb{R}$ of $T$ belong to $\mathcal{C}^1(\overline \Omega)$ and they are not necessarily constant. We also consider the  heat equation  for nonhomogeneous materials with Robin boundary conditions, that consists in finding $v: \Omega\times (0,\tau]\to \rr$ (for $\tau>0$) such that
\begin{equation}\label{PROBgradrobin}
\begin{cases}
&
\partial_t v(\underline{x},t)  +{\rm div}\, T(x) v(\underline{x},t) = 0,\ \ \ \ (\underline{x},t)\in \Omega\times (0,\tau]
\\
&
v(\underline{x},0)=f(\underline{x}),\ \ \ \underline{x}\in \Omega
\\
&
b(\underline x) v(\underline x,t)+ \sum_{\ell=1}^3a_\ell(\underline x)n_\ell(\underline x) \partial_{x_\ell}v(\underline x,t)=0,\ \ (\underline{x},t)\in \partial\Omega\times (0,\tau],
\end{cases}
\end{equation}
where
$n=(n_1,n_2,n_3)$
is the outward unit normal vector to $\partial\Omega$,
and $b:\partial\Omega \to \mathbb{R}$ is a given continuous function. From the physical point of view, if we call $q(\underline x,t)$ the flux  of the quantity described by $v(\underline x,t)$ at the instant $t$, the Fourier's law states that $q(\underline x, t)=T(x) (v(x,t))$.

\medskip
The simpler case is when  we consider $\Omega=\mathbb{R}^3$ and the homogeneous
diffusion problem is the consequence of Fourier's law
$$ q(\underline x,t)=-\nabla v(\underline x,t)$$
and of the conservation of the energy
$$ \partial_t v(\underline x,t)+\operatorname{div}(q(\underline x,t))=0.$$
In this case $T$ is reduced to the negative gradient operator  $T=-\nabla$ and
observing that $\operatorname{div}\circ \nabla = \Delta$
the fractional diffusion model is obtained by replacing in the heat equation
the Laplace operator
by its fractional powers
$$
(-\Delta)^\alpha u(\underline x):= \int_{\rr^3}  \frac{u(\underline x)-u(\underline y)}{|\underline x-\underline y|^{3+2\alpha}}\,dV(\underline y), \ \ \ {\rm for}\ \ \ \alpha\in (0,1).
$$
The fractional versions of the evolution equation in $\rr^3$ is given by
\begin{equation}\label{PROB_grad_frac}
\partial_t v(\underline{x},t)  +(-\Delta)^\alpha v(\underline{x},t) = 0.
\end{equation}

We observe that the fractional diffusion problem modifies both the  Fourier's law and the conservation of the energy.
Using the quaternionic functional calculus  we  are able to define the fractional powers of vector operators, such as $\nabla$ or $T$, in a bounded or unbounded domain $\Omega$ of $\rr^3$. Denoted, just for the moment, by $T^\alpha$ or $\nabla^\alpha$ these fractional operators, we can define the fractional diffusion problem \eqref{PROB_grad_frac} in the divergence form
\begin{equation}\label{PROBgraddirichletfracq}
\partial_t v(\underline{x},t)  +\operatorname{div} T^\alpha(x) v(\underline{x},t) = 0,\ \ \ \ (\underline{x},t)\in \Omega\times (0,\tau].
\end{equation}
The boundary conditions one has to associate with the fractional evolution problem are a very delicate issue and will not discussed here.
We just mention that the most natural boundary conditions are $v=0$ at infinity in the case
$\Omega=\mathbb{R}^3$.

\begin{remark} The boundary condition that
 we have to assume to generate the fractional powers of $T$ are given by
\begin{equation}\label{PROBgradrobinfracq}
a(\underline x) v(\underline x,t)+ \sum_{\ell=1}^3a^2_\ell(\underline x)n_\ell(\underline x) \partial_{x_\ell}v(\underline x,t)=0,\ \ (\underline{x},t)\in \partial\Omega\times (0,\tau],
\end{equation}
where $a:\partial\Omega \to \mathbb{R}$ is a given continuous function.
These Robin-like boundary condition
 differs from the boundary condition in \eqref{PROBgradrobin} by a power two on the coefficients $a_\ell$'s. This is due to the fact that in \eqref{PROBgradrobin} the boundary condition rise from a physical condition on the flux through the boundary, instead the boundary condition in \eqref{PROBgradrobinfracq} naturally comes from the definition of $T^\alpha$. In any case, the two boundary conditions are related when $T$ has coefficients that become constant on the boundary $\partial\Omega$. Indeed, suppose that there exists a constant $\mu$ such that the functions $a_1$, $a_2$, $a_3$  satisfy the conditions
 \begin{equation}\label{VBN}
 a_1(\underline x)=a_2(\underline x)=a_3(\underline x)=\mu \ \ {\rm for\ all}\ \ x\in \partial\Omega
 \end{equation}
  and the coefficients $a$ and $b$ are such that
  \begin{equation}\label{VBNNN}
 a(\underline x)=\mu b(\underline x) \ \ {\rm for\ all}\ \ \underline x\in \partial\Omega.
 \end{equation}
  Then the relation
  $$\sum_{\ell=1}^3a_\ell(\underline x) n_\ell(\underline x) \partial_{x_\ell} v(\underline x) +b(\underline x)v(\underline x)=0$$
is equivalent to
$$\sum_{\ell=1}^3a^2_\ell(\underline x)n_\ell(\underline x) \partial_{x_\ell} v(\underline x) +a(x)v(\underline x)=0$$
when $x\in \partial\Omega.$ For, using (\ref{VBN}) and (\ref{VBNNN}), we have
\begin{equation}\nonumber
\begin{split}
\sum_{\ell=1}^3a^2_\ell(\underline x)n_\ell(\underline x) \partial_{x_\ell}+a(\underline x)I
&
=\mu^2\sum_{\ell=1}^3n_\ell(\underline x) \partial_{x_\ell}+\mu b(\underline x)I
=\mu \Big(\sum_{\ell=1}^3a_\ell(\underline x) n_\ell(\underline x) \partial_{x_\ell}+b(\underline x)I\Big).
\end{split}
\end{equation}
\end{remark}
This kind of approach has several advantages.
\begin{itemize}
\item It generates the fractional Fourier law from the Fourier law
\begin{equation}\label{FFL}
 q(\underline x,t)=T^\alpha(\underline x)v(x,t)
 \end{equation}
using the boundary conditions of the problem and
without modifying the conservation of energy law.
\item
We can define the fractional heat equation for nonhomogeneous materials.\\
\item
The fractional differential equation remains in the divergence form, so the definition of a weak solution in obtained in a simple way.
\item
    It turns out that the approach through the quaternionic functional calculus for defining the fractional heat equations is consistent with the classical one. Indeed, we have that for any $\alpha\in (0,1)$
$$ 2\operatorname{div}(\nabla^\alpha)=(-\Delta)^{\frac 12+\frac\alpha 2}. $$
\end{itemize}

\medskip
Now we present how the quaternions can be used to describe the vector operators. Let $e_\ell$, for $\ell=1,2,3$, be an orthogonal basis for the quaternions $\mathbb{H}$. We identify the vector operator $T$, described in \eqref{T}, with the quaternionic gradient operator with non constant coefficients
\begin{equation}\label{TQ}
\left(
\begin{split}
a_1(\underline x)\partial_{x_1}\\
a_2(\underline x)\partial_{x_2}\\
a_3(\underline x)\partial_{x_3}
\end{split}
\right)
\equiv\sum_{\ell=1}^3e_\ell T_\ell,
\end{equation}
where the components $T_\ell$, $\ell=1,2,3$,  are defined by
$T_\ell:=a_\ell(\underline x)\partial_{x_\ell}$, $\underline x\in \overline{\Omega}$. From the physical point of view the operator $T$, defined in
(\ref{TQ}), can represent the Fourier law for nonhomogeneous materials,
 but it can represent also different physical laws.
Our goal is to generate
 the fractional powers  of $T$, that we denote with $P_\alpha(T)$ for $\alpha\in (0,1)$, when the operators
$T_\ell$, for $\ell=1,2,3$ do not commute among themselves.
\begin{remark}
The notation $P_\alpha(T)$,  for the fractional powers of $T$ ,is more precise with respect to the formal notation $T^\alpha$ is a sense that will be clear just in the following with the precise definition.
 The formal notation $T^\alpha$ is used in \eqref{PROBgraddirichletfracq} and \eqref{PROBgradrobinfracq} (see \eqref{FFL}).
\end{remark}
\begin{definition}
The vector part of the  fractional powers $P_{\alpha}(T)$ is called the fractional Fourier law associated with $T$.
\end{definition}
  Now we present the general theory of the $S$-spectrum to construct the fractional power of a quaternionic right linear operator.
Let $V$ be a two-sided quaternionic Banach space and $\closOP(V)$ the set of closed quaternionic right linear operators on $V$. The Banach space of all bounded right linear operators on $V$ is indicated by the symbol $\mathcal{B}(V)$ and is endowed with the natural operator norm.
For $T\in\closOP(V)$, we define the operator associated with the $S$-spectrum as:
\begin{equation}\label{QST}
\Q_{s}(T) := T^2 - 2\Re(s)T + |s|^2\id, \qquad \text{for $s\in\hh$}
\end{equation}
where $\Q_{s}(T):dom(T^2)\to V$, where $dom(T^2)$ is the domain of $T^2$.
 We define the $S$-resolvent set of  $T$ as
\[\rho_S(T):= \{ s\in\hh: \Q_{s}(T) \ {\rm is\ invertible\ and \ } \Q_{s}(T)^{-1}\in\boundOP(V)\}\]
and the $S$-spectrum of $T$ as
\[\sigma_S(T):=\hh\setminus\rho_S(T).\]
The operator $\Q_{s}(T)^{-1}$ is called the pseudo $S$-resolvent operator.
 For $s\in\rho_S(T)$, the left $S$-resolvent operator is defined as
\begin{equation}\label{SRESL}
S_L^{-1}(s,T):= \Q_s(T)^{-1}\overline{s} -T\Q_s(T)^{-1}
\end{equation}
and the right $S$-resolvent operator is given by
\begin{equation}\label{SRESR}
S_R^{-1}(s,T):=-(T-\id \overline{s})\Q_s(T)^{-1}.
\end{equation}
The fractional powers of $T$, denoted by $P_{\alpha}(T)$, are defined as follows:
for any $I\in \mathbb{S}$,  for $\alpha\in(0,1)$ and  $v\in\dom(T)$ we set
\begin{equation}\label{BALA1}
P_{\alpha}(T)v := \frac{1}{2\pi} \int_{-I\rr}   S_L^{-1}(s,T)\,ds_I\, s^{\alpha-1} T v,
\end{equation}
or
\begin{equation}\label{BALA2}
P_{\alpha}(T)v := \frac{1}{2\pi} \int_{-I\rr} s^{\alpha-1} \,ds_I\,  S_R^{-1}(s,T) T v,
\end{equation}
where $ds_j=ds/I$.
These formulas are a consequence
of the quaternionic version of the $H^\infty$-functional calculus based on the $S$-spectrum, see the book \cite{6CG} for more details.
For the generation of the fractional powers $P_{\alpha}(T)$ a crucial assumption
on the $S$-resolvent operators is that, for $s\in \mathbb{H}\setminus \{0\}$ with ${\rm Re}(s)=0$, the estimates
\begin{equation}\label{SREST}
\left\|S_L^{-1}(s,T)\right\|_{\mathcal{B}(V)} \leq \frac{\Theta}{|s|}\quad\text{and}\quad
\left\|S_R^{-1}(s,T)\right\|_{\mathcal{B}(V)} \leq \frac{\Theta}{|s|},
\end{equation}
hold
with a constant $\Theta >0$ that does not depend on the quaternion $s$.
 It is important to observe that
the conditions (\ref{SREST}) assure that the integrals (\ref{BALA1}) and (\ref{BALA2}) are convergent and so the fractional powers are well defined.

For the definition of the fractional powers of the operator $T$ we can use equivalently the integral representation in (\ref{BALA1}) or the one in (\ref{BALA2}).
Moreover, they correspond to a modified version of Balakrishnan's formula that takes only spectral points with positive real part into account.

\medskip
We want to apply the previous theory to the case $V:= L^2(\Omega, \hh)$ and $T\in\mathcal K(L^2(\Omega,\hh))$ defined as in \eqref{TQ} ($\operatorname{dom}(T)\subset L^2(\Omega, \hh)$ is a densely subset). A crucial problem is to determine the conditions on the coefficients $a_1$, $a_2$, $a_3:\overline{\Omega} \subset\mathbb{R}^3\to \mathbb{R}$ such that \eqref{BALA1} and \eqref{BALA2} are convergent. This problem is splitted into two problems
\begin{itemize}
\item the first is to find appropriate conditions for the coefficients $a_i$'s such that the purely imaginary quaternions are in the $S$-resolvent set $\rho_S(T)$ (i.e. $\mathcal Q_s(T): \operatorname{dom}(T^2)\to L^2(\Omega,\hh)$ is invertible and bounded).
This is a necessary condition, see formulas (\ref{BALA1}) and (\ref{BALA2}). Then, since
in the quaternionic case the map
$s\mapsto s^\alpha$, for $\alpha\in (0,1)$ is not defined for $s\in(-\infty,0)$ and, unlike in the complex setting, it is not possible to choose different branches of $s^{\alpha}$ in order to avoid this problem.
For this reason it is of great importance to assume the condition $\Re(s) \geq 0$  that avoids the
half real line $(-\infty,0]$.
\\
\item The second crucial fact is to determine the conditions on the coefficients $a_i$'s such that the estimate \eqref{SREST} for the $\mathcal S$-resolvent operator of $T$ holds true.
\end{itemize}
Both these problems are solved by considering the following approach. According to the initial condition of the  boundary-value problems we invert the operator $\mathcal Q_s(T)$ on the space $H^1_0(\Omega,\hh)$, when we consider the Dirichlet boundary condition, and on the space
$$
\mathcal H:=\{ u\in H^1(\Omega)|\, \int_\Omega u(\underline x)\, dV(\underline x)=0\, \}
$$
 when we consider the Robin boundary condition. The invertibility of $\mathcal Q_s(T)$ is thus reduced to solve in a weak sense the following two partial differential equations: given $F\in L^2(\Omega,\hh)$ and $s\in \hh \setminus\{0\}$ such that $\operatorname{Re}(s)=0$
\begin{equation}\label{Prob1}
\begin{cases}
 &Q_s(T)(u)=\big( T^2  -2s_0T+ |s|^2\id\big)u(\underline x)=F(\underline x),\ \ \ \underline x\in \Omega,
 \\
 &
u\in H^1_0(\Omega,\hh),
\end{cases}
\end{equation}
and
\begin{equation}\label{Prob2}
\begin{cases}
 &Q_s(T)(u)=\big( T^2  -2s_0T+ |s|^2\id\big)u(\underline x)=F(\underline x),\ \ \ \underline x\in \Omega,
 \\
 &
u\in\mathcal H(\Omega,\hh),
\\
&
b(\underline x) v(\underline x)+ \sum_{\ell=1}^3a^2_\ell(\underline x)n_\ell(\underline x) \partial_{x_\ell}v(\underline x)=0,\ \ \underline x\in \partial \Omega.
\end{cases}
\end{equation}
To solve in the weak sense \eqref{Prob1} (resp. \eqref{Prob2}) means that for any $F\in L^2(\Omega,\hh)$
we have to find $u_F\in H^1_0(\Omega,\hh)$ (resp. $u_F\in\mathcal H(\Omega,\hh)$) such that:
 for any $v\in H^1_0(\Omega,\hh)$ (resp. $v\in\mathcal H(\Omega,\hh)$) we have
\begin{equation}\label{weak}
\langle Q_s(T)(u_F),v \rangle=(F,v)_{L^2}=\int_\Omega \overline F v\, dV(\underline x),
\end{equation}
where the angle-brackets means that $\mathcal Q_s(T)$ is applied to $u_F$ in the sense of distribution. If we solve \eqref{weak}, we can define
$$ Q_s(T)^{-1}(F):= u_F $$
In order to solve \eqref{weak} we  apply the Lax-Milgram Lemma to the sesquilinear form: $\langle Q_s(T)(u_F),v \rangle$. Thus it is crucial to prove an explicit formula for the left hand side of \eqref{weak}. This formula can be deduced from an arguments of integration by parts and using the Dirichlet boundary condition for \eqref{Prob1}:
\begin{equation}\label{b1}
\begin{split}
& \langle Q_s(T)(u),v \rangle=  \sum_{\ell = 1}^{3}\int_{\Omega}   \overline{a_\ell(\underline x)\partial_{x_\ell}(u(\underline x))}\, a_\ell(\underline x)\partial_{x_\ell}(v(\underline x))\,dV(\underline x) +
 \\ & \frac 12\sum_{\ell = 1}^{3} \int_{\Omega}    \overline{\partial_{x_\ell}(u(\underline x))} \partial_{x_{\ell}} \left(a^2_{\ell}(\underline x)\right)v(\underline x)\,dV(\underline x)
+ ( {\rm Vect}(Q_{s}(T)) u,v )_{L^2}
\\ &
+ |s|^2 ( u,v)_{L^2}(:=b_{s,1}(u,v))
\end{split}
\end{equation}
or the Robin boundary condition for \eqref{Prob2}
\begin{equation}\label{b2}
\begin{split}
 & \langle Q_s(T)(u),v \rangle =  \sum_{\ell = 1}^{3}\int_{\Omega}   \overline{a_\ell(\underline x)\partial_{x_\ell}(u(\underline x))}\, a_\ell(\underline x)\partial_{x_\ell}(v(\underline x))\,dV(\underline x)\\
 & +\frac 12\sum_{\ell = 1}^{3} \int_{\Omega}    \overline{\partial_{x_\ell}(u(\underline x))} \partial_{x_{\ell}} \left(a^2_{\ell}(\underline x)\right)v(\underline x)\,dV(\underline x)
+({\rm Vect}(Q_{s}(T)) u,v )_{L^2}
\\
&
+\int_{\partial \Omega} a(\underline x) \overline{u(\underline x)} v(\underline x)\,dS(\underline x)+ |s|^2 ( u,v)_{L^2}
(=:b_{s,2}(u,v)).
\end{split}
\end{equation}

    In \cite{frac1}, \cite{frac2}, \cite{frac3}, \cite{frac4} and \cite{frac5} we found suitable conditions for the coefficients $a_i$'s such that  the two sesquilinear forms $b_{s,1}$ and $b_{s,2}$ are coercive and continuous when $\Omega$ is bounded or unbounded. In conclusion by the Lax-Milgram lemma we obtain the solvability of the equation \eqref{weak} (i.e. the invertibility of $\mathcal Q_s(t)$) and the estimate \eqref{SREST} for the $\mathcal S$-resolvent operator. Thus the problem of the convergence of \eqref{BALA1} and \eqref{BALA2} is solved. In the next three paragraphs we summarize the conditions we found on the coefficients of $T$ to obtain the convergence of \eqref{BALA1} and \eqref{BALA2} .

\medskip
{\em The fractional Fourier's law in the problem \eqref{PROBgraddirichletfracq} with $\Omega$ bounded.}
In the paper \cite{frac1} it was considered the commutative Fourier's law $T_{com}$, that is an operator of the form
\begin{equation}\label{TCOM}
T_{com}=a_1(x_1)\partial_{x_1}e_1 + a_2(x_2)\partial_{x_2}e_2 + a_3(x_3)\partial_{x_3}e_3
\end{equation}
where the real operators
$a_1(x)\partial_{x_1}$, $a_2(x_2)\partial_{x_2}$ and $a_3(x_3)\partial_{x_3}$ commute among themselves.
It has been shown that if the coefficients  $a_{\ell}:\overline{\Omega} \to \mathbb{R}$, for $\ell = 1,2,3$ belong to
$ \mathcal{C}^1(\overline{\Omega},\rr)$ and
if $a_{\ell}$, for $\ell = 1,2,3$ are suitably large and the
 their derivative are suitably small then the integrals in \eqref{BALA1} and \eqref{BALA2} are convergent.

In \cite{frac2} we replace the commutative Fourier's law
$T_{com}$ by the more general Fourier's law
\begin{equation}\label{tncom}
T(x) = a_1(\underline{x})\partial_{x_1} e_1 + a_2(\underline{x}) \partial_{x_2}e_2 + a_3(\underline{x})\partial_{x_3}e_3
\end{equation}
where now the real operators $a_1(\underline{x})\partial_{x_1}$, $ a_2(\underline{x}) \partial_{x_2}e_2$ and $a_3(\underline{x})\partial_{x_3}$ do not commute among themselves. In this case the conditions for the existence of the fractional powers are more complicated.

The main result is summarized in the following theorem (see for more details
Theorems $4.1$, $4.4$ and $4.5$ in \cite{frac2}).

\medskip

\begin{theorem}
Let $\Omega$ be a bounded $\mathcal{C}^1$-domain in $\mathbb{R}^3$,
let $T=\sum_{i=1}^3 a_i(\underline{x}) \partial_{x_i} e_i$
with  $a_i\in \mathcal{C}^1(\overline{\Omega})$ for any $i=1,\, 2,\, 3$
and set
$$
F_{(a_1,a_2,a_3)} :=\sum_{i=1}^3 e_i\partial_{x_i}(a_i ).
$$
Let $a_1,a_2,a_3 \geq m > 0$,  and assume that
\begin{equation}
\min\{\inf_{\underline{x}\in \Omega} a_1^2,\inf_{\underline{x}\in \Omega}a_2^2,\inf_{\underline{x}\in \Omega}a_3^2\}
-(2\max\{\sup_{\underline{x}\in \Omega} a_1^2,\sup_{\underline{x}\in \Omega} a_2^2,\sup_{\underline{x}\in \Omega} a_3^2\})^{1/2}C_\Omega\|F_{(a_1,a_2,a_3)}\|_{L^\infty}>0
\end{equation}
and
\begin{equation}
1-2 \|F_{(a_1,a_2,a_3)}\|_{L^\infty}\Big(1+ 4C^2_\Omega
\max\{\sup_{\underline{x}\in \Omega}(1/a_1^2),\sup_{\underline{x}\in \Omega}(1/a_2^2),\sup_{\underline{x}\in \Omega}(1/a_3^2)\}
\Big)>0,
\end{equation}
where $C_\Omega$ in the Poincar\'{e} constant of $\Omega$ and
$$
 \|F_{(a_1,a_2,a_3)}\|_{L^\infty}:=\sup_{\underline{x}\in \Omega}( |\partial_{x_1}(a_1)|+|\partial_{x_2}(a_2)|+|\partial_{x_3}(a_3)|).
 $$
Then for any $\alpha\in(0,1)$ and for any  $v\in\dom(T)$, the integrals \eqref{BALA1} and \eqref{BALA2}
converge absolutely.
\end{theorem}

\medskip
The sesquilinear form $b_{s,1}(u,v)$, defined in \eqref{b1},
associated with the invertibility of the operator
$\Q_s(T):=T^2-2s_0T+|s|^2\mathcal{I}$, with homogeneous Dirichlet
boundary conditions, has to be considered with care.
 We summarized in Remark \ref{RKfacts} some considerations associated with $b_{s,1}(u,v)$.

 \begin{remark}\label{RKfacts}
 We point out some fact that appear in the application of the Lax--Milgram lemma according to the dimension $n$
  of $\Omega$.

 (I) The quadratic form $b_{s,1}(u,v)$
  associated to the operator $\mathcal{Q}_s(T)$  is in general degenerate on $H^1_0(\Omega,\mathbb{H})$.

(II)
In dimension $n=3$ , when $\Omega$ is a $\mathcal{C}^1$ bounded set in $\mathbb{R}^3$ and
 $a_1\not=0$, $a_2\not=0$, $a_3\not=0$, it turns out that $b_{s,1}(u,v)$ is continuous and coercive
under suitable conditions on the coefficients
of $\mathcal{Y}\times \mathcal{Y}$,
where
$$
 \mathcal{Y}:=\{v\in H^1_0(\Omega,\mathbb{H}) \ : \ v_0=v_1=v_2=v_3\}
 $$
 is a closed subspace of $H^1_0(\Omega,\mathbb{H})$ and the
 $S$-resolvent operators satisfy suitable growth conditions which
 ensure the existence of the fractional powers.

(III)
In dimension $n=2$, when $\Omega$ is a $\mathcal{C}^1$ bounded set in $\mathbb{R}^2$
and $a_1\not=0$, $a_2\not=0$, it turns out that $b_{s,1}(u,v)$ is continuous and coercive
under suitable conditions on the coefficients
of $\mathcal{X}\times \mathcal{X}$,
where
$$
\mathcal{X}:=\{v\in H^1_0(\Omega,\mathbb{H})\ :\ \ v_0=v_2\ \ {\rm and } \  \ v_1=v_3 \}
$$
 is a closed subspace of $H^1_0(\Omega,\mathbb{H})$ and the
 $S$-resolvent operators satisfy suitable growth conditions which
 ensure the existence of the fractional powers.

(IV)
If we consider the quadratic form
in dimension $n=3$, that is when $\Omega$ is a $\mathcal{C}^1$ bounded set in $\mathbb{R}^3$
and $a_1\not=0$, $a_2\not=0$, $a_3=0$, then the quadratic form is not coercive because of $a_3=0$.
It seems that this case cannot be treated using Lax--Milgram Lemma, but a suitable method for degenerate equations has to be used.

(V)
In \cite{frac2} we proved, under more restrictive hypothesis on the coefficients $a_i$'s, that the sesquilinear form $b_{s,1}(u,v)$ is continuous and coercive in $H^1_0(\Omega,\hh)$.

(VI)
From the physical point of view the case for $a_1\not=0$, $a_2\not=0$, $a_3=0$, in dimension $n=3$ is the
case in which the conductivity is the direction $z$ goes to zero.

(VII)
The proofs for the continuity and coercivity are similar in any dimension,
the estimate for the $S$-resolvent operators have some differences
according to the fact that we work in $\mathcal{Y}$ or in $\mathcal{X}$.
\end{remark}

\medskip
{\em The fractional Fourier's law in the problem \eqref{PROBgradrobinfracq} with $\Omega$ bounded.}
Regarding the initial boundary value problem of Robin-type, the conditions on the coefficients $a_i$'s depend on two positive constants which appear in the following two inequalities:
\begin{itemize}
\item for all $u\in H^1(\Omega, \mathbb{R})$ the following inequality holds:
\begin{equation}\label{trace}
\|u\|_{H^{1/2}(\partial\Omega,\mathbb{R})}\leq C_{\partial\Omega}\|u\|_{H^1(\Omega,\mathbb{R})}.
\end{equation}
where $C_{\partial\Omega}$ does not depend on $u$
\item for all $u\in H^1(\Omega, \mathbb{R})$ the following inequality holds:
\begin{equation}\label{PWI}
\left\|u-|\Omega|^{-1}\int_\Omega u(x)dx\right\|_{L^2(\Omega,\mathbb{R})}\leq C_P\|\nabla u\|_{L^2(\Omega,\mathbb{R})},
\end{equation}
where $C_P$ does not depend on $u$.
\end{itemize}
We obtained in \cite{frac5} the following result (see Theorems $4.4$, $5.1$ and $5.2$ in \cite{frac5}).
\begin{theorem}
Let $\Omega$ be a bounded domain in $\mathbb R^3$ with boundary $\partial\Omega$ of class $\mathcal C^1$. Assume that
$a\in \mathcal{C}^0(\partial\Omega, \mathbb{R}) $ and
 let $T$ be the operator defined in (\ref{T}) with coefficients $a_1$, $a_2$, $a_3\in \mathcal{C}^1(\overline{\Omega}, \mathbb{R})$.
 Define the following constants:
\begin{equation}\nonumber
C_T:=\min_{\ell=1,2,3}\inf_{ x\in\Omega}(a^2_\ell(x)),\quad
C_T':=\sum_{i,\ell=1}^3\|a_\ell\partial_{x_\ell}a_i\|_\infty,\quad K_{a,\, \Omega}:= C^2_{\partial\Omega}\|a\|_\infty,
\end{equation}
where  $\| \cdot\|_\infty$ denotes the sup norm and $ C_{\partial\Omega}$ is the constant in \eqref{trace}. Moreover, assume that
\begin{equation}\label{kappaomeg}
 C_T-C_T' C_P-K_{a, \Omega}\Big(1+C_P^2 \Big)>0 \ \ \ and \ \ \ \ \
C_T>0,
\end{equation}
where $C_P$ is the constant in (\ref{PWI}).
Then for any $\alpha\in(0,1)$ and for any  $v\in\dom(T)$, the integrals \eqref{BALA1} and \eqref{BALA2}
converge absolutely.
\end{theorem}
\begin{remark}
In \cite{frac4} we treated the case of the operator $T$ with commutative coefficients and of $\Omega$ bounded with Robin-type boundary condition.
\end{remark}

{\em The fractional Fourier's law in the problem \eqref{PROBgraddirichletfracq} with $\Omega$ unbounded.} Regarding the initial boundary value problem of Dirichlet-type for the unbounded domains, we obtained in \cite{frac5} the following result (see Theorems $4.8$, $5.1$ and $5.2$ in \cite{frac5}).
\begin{theorem}
Let $\Omega$ be an unbounded domain in $\mathbb R^3$ with boundary $\partial\Omega$ of class $\mathcal C^1$. Let $T$ be the operator defined in \eqref{T} with coefficients $a_1$, $a_2$, $a_3\in \mathcal{C}^1(\overline{\Omega}, \mathbb{R})$. Suppose that
\begin{equation}\label{c1bis}
M:=\sum_{i,j=1}^3\|a_i\partial_{x_i}(a_j)\|_{L^3(\Omega)}< +\infty
\end{equation}
and
\begin{equation}\label{kappaomegbis}
 C_T:=\min_{\ell=1,2,3}\inf_{\underline x\in\Omega} (a^2_\ell(\underline x))>0, \ \ \ \
  C_T-4M>0.
\end{equation}
Then for any $\alpha\in(0,1)$ and for any  $v\in\dom(T)$, the integrals \eqref{BALA1} and \eqref{BALA2}
converge absolutely.
\end{theorem}

\end{document}